\documentclass{amsart}
\usepackage{amsmath,amssymb,amsfonts}
\usepackage{diagrams}
\diagramstyle[nohug]

\newtheorem{thm}{Theorem}[section]
\newtheorem{cor}[thm]{Corollary}
\newtheorem{example}[thm]{Example}
\newtheorem{lem}[thm]{Lemma}
\newtheorem{prop}[thm]{Proposition}
\newtheorem{defn}[thm]{Definition}
\newtheorem*{linear}{Linearity Law}
\newtheorem*{vanish}{Vanishing Law}
\newtheorem*{transform}{Transformation Law}

\def\gfrac#1#2{\left[\begin{array}{c}#1\\#2\end{array}\right]}

\newcommand{\fP}{{\mathfrak P}}
\newcommand{\fQ}{{\mathfrak Q}}
\newcommand{\m}{{\mathfrak m}}
\newcommand{\p}{{\mathfrak p}}
\newcommand{\q}{{\mathfrak q}}
\newcommand{\M}{{\mathcal M}}
\DeclareMathOperator{\height}{ht}
\DeclareMathOperator{\mH}{H}
\DeclareMathOperator{\Hom}{Hom}
\DeclareMathOperator{\Ext}{Ext}
\DeclareMathOperator{\Spec}{Spec}

\begin{document}


\title{Module structure of an injective resolution}
\author{C-Y. Jean Chan \and I-Chiau Huang}
\address{Department of Mathematical Sciences, University of Arkansas,
                 Fayetteville, AR 72701, U.S.A. \and
                 Institute of Mathematics, Academia Sinica, Nankang,
                 Taipei 11529, Taiwan, R.O.C.}
\email{cchan@uark.edu \and ichuang@math.sinica.edu.tw}

\begin{abstract}
   Let $A$ be the ring obtained by localizing the polynomial ring
   $\kappa[X,Y,Z,W]$ over a field $\kappa$
   at the maximal ideal $(X,Y,Z,W)$ and modulo the
   ideal $(XW-YZ)$. Let $\p$ be the ideal of $A$ generated by $X$ and
   $Y$.  We study the module structure of a minimal injective
   resolution of $A/\p$ in details using local cohomology.
   Applications include the description of $\Ext^i_A(M, A/\p)$, where $M$
   is a module constructed by Dutta, Hochster and
   McLaughlin, and the Yoneda product of
   $\Ext^*_A(A/\p, A/\p)$.
\end{abstract}
\maketitle


\section{introduction}\label{sec:intro}


In the category of modules over a commutative ring, injective and
projective modules are dual notions.  To study cohomology properties
of a module, we may consider a minimal free resolution or a minimal
injective resolution of the module. The boundary maps of the former
are given by matrices in terms of given basis. The coboundary maps of
the latter are discussed less extensively. In general, there are no
simple descriptions for injective resolutions.  The subtlety comes
partly from the fact that there are no canonical ways to identify
minimal injective modules (injective hulls) for a given module, even
though they are all isomorphic. On works regarding concrete realizations
of Grothendieck duality, one finds many natural injective hulls with different 
guises for a given module. The difference of these injective hulls is a part of the 
structure of the underlying module. From this viewpoint,
injective hulls for a given module are not unique, just as there are
different vector spaces of the same dimension (for instance, a finite
dimensional vector space and its dual). Intriguing structures such as
residues support this viewpoint since they arise from
isomorphisms between injective hulls.

Here is a typical example: Let $R$ be a formal power series ring of
$n$ variables over a field $\kappa$. The $n$th local cohomology module
of $R$ supported at the maximal ideal gives rise to an injective hull
of $\kappa$. The elements of this local cohomology module have a 
concrete description using generalized fractions. As an $R$-module, 
$\kappa$ has another injective hull consisting of the $\kappa$-linear 
homomorphisms from $R$ to $\kappa$ annihilated by some power of the 
maximal ideal. Residues appear when an explicit isomorphism
 between these two injective hulls is constructed. The reader is referred to
\cite{hu:pmzds,hu:ecrc} for more details and further developments
along this direction.

The goal of our work is to develop a concrete means to study the
structure of injective resolutions of modules. It consists of two
steps: constructing injective modules explicitly and then describing
the coboundary maps explicitly in a resolution built up from the
injective modules obtained in the previous step.  The goal has been 
achieved for modules related to residual complexes, which are of particular 
interests due to their central role in Grothendieck duality theory. We recall
that residual complexes are build by injective hulls of the residue fields of 
points on a scheme and resolve canonical modules in certain Cohen-Macaulay cases.
In \cite{hu:ecrc}, residual complexes are constructed concretely in a relatively canonical
way.  The construction in
\cite{hu:ecrc} is local. One of its globalizations gives rise to
injective resolutions for the vector bundles on projective spaces
\cite{hu:cpssrc}. The injective resolutions obtained in \cite{hu:ecrc,
hu:cpssrc} are for modules ({\em resp.} sheaves of modules) whose
structures ({\em resp.} local structures) are determined completely by the
underlying rings ({\em resp.} schemes). Not much is known in general
about concrete constructions of injective resolutions of non-flat modules.

Studies of the homology and cohomology modules from the viewpoint of
injective objects are often restricted to some subcategories of the
category of modules, such as the category of graded (or multi-graded)
modules (see for example \cite{got-wat:gr,mil:Adfldms,mil-sturm:cca})
or the category of squarefree modules \cite{yan:smlcmmi}. Injective
resolutions in these smaller categories drastically differ from those
in the category of all modules. For instance, a multi-graded injective
resolution of the polynomial ring $\kappa[X,Y]$ of two variables over
a field $\kappa$ consists of only four indecomposable multi-graded
injective modules \cite[Example~11.20]{mil-sturm:cca}. In the category
of modules concerning no gradings, a minimal injective resolution of
$\kappa[X,Y]$ consists of infinitely many indecomposable injective
modules indexed by the prime ideals of $\kappa[X,Y]$ due to its Gorenstein property
(\cite[Theorem~18.8]{mats:crt}).
Minimal injective resolutions, especially those for modules over a
local ring, are still full
of mysteries and are not
possible to be deduced from graded cases. At the time when more case
studies of minimal injective resolutions are available, a general
theory may be developed for a larger class of modules. This paper
serves as a first step towards such direction by carrying out the
above goal for a module related to an important example in the discussions of several
homological conjectures~({\em c.f.}  \cite{dut-hoc-mcl:mfpdnim}).

In this paper, our study emphasizes the module structure of injective resolutions
rather than its category structure. More precisely, we would like to construct
explicitly an injective resolution of a given module and obtain its
cohomological information from the resolution. Let $S$ be the polynomial ring
$\kappa[X,Y,Z,W]$  over a field $\kappa$. In this paper, we consider the ring
$$
A = S_{(X,Y,Z,W)}/(XW-YZ)
$$
and the ideal $\p$ of $A$ generated by $X,Y$. For each prime ideal
$\q$ of $\kappa[Z,W]$ contained in $\m:=(Z,W)$, we construct an
injective hull $E(A/(\q,X,Y))$ of $A/(\q,X,Y)$. In the sequel, we write
$E(A/(\q,X,Y))$ simply as $E(\q)$. In terms of {\em generalized
   fractions} (defined in Definition~\ref{defn:21}) of elements of
$E(\q)$, our main result describes a minimal injective resolution
$$
E(0)\to\underset{\scriptscriptstyle\q\neq\m}{\oplus}E(\q)
\to\underset{\scriptscriptstyle\q\neq(0)}{\oplus}E(\q)\to E(\m)^2\to 
E(\m)^2\to E(\m)^2\to\cdots
$$
of $A/\p$. According to the authors' knowledge, this is the first detailed
analysis of an injective resolution  for a module which does not come 
from duality theory.

As applications, we read explicitly
\begin{itemize}
\item local cohomology modules $\mH^i_I(A/\p)$ of $A/\p$ supported at 
an ideal $I$ of $A$,
\item an isomorphism $\Hom_A(\p/\p^2,A/\p)\to\Ext^1_A(A/\p,A/\p)$ of 
normal modules,
\item the product of the Yoneda algebra $\Ext^*_A(A/\p,A/\p)$,
\item $\Ext^i_A(M,A/\p)$, where $M$ is the $A$-module constructed by
        Dutta, Hochster and McLaughlin \cite{dut-hoc-mcl:mfpdnim}.
\end{itemize}

The paper is organized as follows: In Section~\ref{sec:genfrac}, we
recall the notion of generalized fractions which describe elements in
certain top local cohomology modules. Technical properties are prepared
for latter use.  In Section~\ref{sec:injmod}, we construct injective
hulls in terms of generalized fractions. These injective
hulls are building blocks for our injective resolution. In
Section~\ref{sec:injresl}, we define homomorphisms for these injective
hulls and show that they give rise to a minimal injective resolution.
In Section~\ref{sec:cohomology}, we carry out the computations
for the applications listed in the previous paragraph.

The authors thank L. Avramov for pointing out an error regarding the description of Yoneda
algebra in an earlier version of this paper.


\section{generalized fractions}\label{sec:genfrac}


Our description of injective modules and coboundary maps of an injective 
resolution is
based on local cohomology modules and the representation of their 
elements by generalized
fractions. We recall the definition and some properties of 
generalized fractions and
refer the details to \cite[Chapter 2]{hu:pmzds}. Let $R$ be a 
Noetherian ring and $I$ be
an ideal of $R$ generated {\em up to radical} by $n$ elements
$x_1,\cdots,x_n$ and another ideal $J$. Let $N$ be an $R$-module, 
whose elements
are annihilated by a power (depending on the element) of $J$. 
Elements of the $n$-th
local cohomology module $\mH^n_I(N)$ of $N$ supported at $I$ can be 
described by
the following exact sequence
\begin{equation}\label{eq:34123}
\oplus_{i=1}^n N_{x_1\cdots\widehat{x_i}\cdots 
x_n}\xrightarrow{\alpha}N_{x_1\cdots x_n}\xrightarrow{\beta}
\mH^n_I(N)\to 0,
\end{equation}
where $\alpha$ is the map given by
$$
\frac{\omega}{(x_1\cdots\widehat{x_i}\cdots x_n)^s}\mapsto
\frac{(-1)^i x_i^s\omega}{(x_1\cdots x_n)^s}
$$
for $\omega\in N$ and $s\geq 0$.
\begin{defn}\label{defn:21}
A generalized fraction
$$
\gfrac{\omega}{x_1^{i_1},\cdots,x_n^{i_n}}\in\mH^n_I(N),
$$
where $\omega\in N$ and $i_1,\cdots,i_n\in\mathbb Z$, is the image of
$x_1^{s-i_1}\cdots x_n^{s-i_n}\omega/(x_1\cdots x_n)^s$ under the
map $\beta$ in (\ref{eq:34123}) for a sufficiently large $s$. 
$\omega$ is called the
numerator of the generalized fraction and 
$x_1^{i_1},\cdots,x_n^{i_n}$ are called the
denominators of the generalized fraction.
\end{defn}
If some $i_j$ is less than one, then the above generalized fraction 
vanishes. Generalized
fractions satisfy the following properties.

\begin{linear}
For $\omega_1,\omega_2 \in N$ and $a_1,a_2 \in R$,
$$
\gfrac{a_1\omega_1+a_2\omega_2}{x_1,\cdots,x_n}
=
a_1\gfrac{\omega_1}{x_1,\cdots,x_n}+a_2\gfrac{\omega_2}{x_1,\cdots,x_n}.
$$
\end{linear}

\begin{transform}
For $\omega \in N$ and elements $x_1', \cdots, x_n'$, which together with $J$
generate $I$ up to radical,
$$
\gfrac{\omega}{x_1,\cdots,x_n}=\gfrac{\det(r_{ij})\omega}{x_1',\cdots,x'_n}
$$
if $x_i' = \sum_{j=1}^{n} r_{ij} x_j$ for $i=1,\cdots,n$.
\end{transform}

\begin{vanish}
For $\omega\in N$,
$$
\gfrac{\omega}{x_1,\cdots,x_n}=0
$$
if and only if $(x_1\cdots x_n)^s\omega\in (x_1^{s+1},\cdots,x_n^{s+1})N$
for some $s\geq 0$.
we can take $s=0$.
\end{vanish}
Note that powers of $x_1,\cdots,x_n$ together with $J$ also generate $I$ up to
radical. So the above laws apply to generalized fractions with 
arbitrary denominators.
An easy application of these laws is that adding to one of the 
denominators by a linear
combination of other denominators does not change the value.
\begin{example}\label{example2.2}
Look at the case $n=2$. $I$ is generated by $x_1,x_2-ax_1$ and $J$ up 
to radical for any $a\in R$. We have
\begin{eqnarray*}
\gfrac{\omega}{x_1,x_2-ax_1}&=&\gfrac{x_2\omega}{x_1,x_2(x_2-ax_1)}\\
&=&\gfrac{(x_2-ax_1)\omega}{x_1,x_2(x_2-ax_1)}+\gfrac{ax_1\omega}{x_1,x_2(x_2-ax_1)}\\
&=&\gfrac{\omega}{x_1,x_2}.
\end{eqnarray*}
\end{example}
\begin{prop}\label{prop:38555}
Let $R$ be a Noetherian local ring and $N$ be an $R$-module, whose elements are
annihilated by a power (depending on the element) of the maximal 
ideal of $R$. An
element of 
$\mH^n_{(X_1,\cdots,X_n)}(R[X_1,\cdots,X_n]_{(X_1,\cdots,X_n)}\otimes 
N)$ can
be written as
\begin{equation}\label{eq:89388}
\Psi=\sum_{i_1,\cdots,i_n>0}
\gfrac{1\otimes\alpha_{i_1\cdots i_n}}{X_1^{i_1},\cdots,X_n^{i_n}},
\end{equation}
where $\alpha_{i_1\cdots i_n}\in N$. The expression is unique in the 
sense that $\Psi=0$ if and
only if $\alpha_{i_1\cdots i_n}=0$ for all $i_1\cdots i_n>0$.
\end{prop}
\begin{proof}
$N$ has a natural module structure over the completion $\hat{R}$ of $R$.
Elements of $\mH^n_{(X_1,\cdots,X_n)}(\hat{R}[[X_1,\cdots,X_n]]\otimes N)$ can
be written uniquely in the form of (\ref{eq:89388}), see \cite[p. 
21]{hu:pmzds}.
The proposition follows from the canonical isomorphism
$$
\mH^n_{(X_1,\cdots,X_n)}(R[X_1,\cdots,X_n]_{(X_1,\cdots,X_n)}\otimes N)
\simeq\mH^n_{(X_1,\cdots,X_n)}(\hat{R}[[X_1,\cdots,X_n]]\otimes N).
$$
\end{proof}
Let $S$ be the polynomial ring  $\kappa[X,Y,Z,W]$ over a field $\kappa$
  as in Section~\ref{sec:intro}.
\begin{cor}\label{cor:ZW}
Elements of $\mH^4_{(X,Y,Z,W)}(S_{(X,Y,Z,W)})$ can be written uniquely as
\begin{equation}\label{eq:967813}
\sum_{i,j,k,\ell>0}\gfrac{a_{ijk\ell}}{Z^i,W^j,X^k,Y^\ell},
\end{equation}
where $a_{ijk\ell}\in\kappa$.
\end{cor}
We call $a_{ijk\ell}$ the coefficient of 
$\gfrac{1}{Z^i,W^j,X^k,Y^\ell}$ for the
element (\ref{eq:967813}).
\begin{cor}\label{cor:0}
Elements of $\mH^2_{(X,Y)}(S_{(X,Y)})$ can be written uniquely as
$$
\sum_{i,j>0}\gfrac{\varphi_{ij}}{(XW)^i,(YZ)^j},
$$
where $\varphi_{ij}\in\kappa(Z,W)$.
\end{cor}
\begin{proof}
$S_{(X,Y)}\simeq\kappa(Z,W)[XW,YZ]_{(XW,YZ)}$.
\end{proof}
\begin{cor}\label{cor:Z}
Elements of $\mH^3_{(X,Y,Z)}(S_{(X,Y,Z)})$ can be written uniquely as
$$
\sum_{i,j,k>0}\gfrac{\varphi_{ijk}}{Z^i,(XW)^j,Y^k},
$$
where $\varphi_{ijk}\in\kappa(W)$.
\end{cor}
\begin{proof}
$S_{(X,Y,Z)}\simeq\kappa(W)[XW,Y,Z]_{(XW,Y,Z)}$.
\end{proof}
\begin{cor}\label{cor:W}
Elements of $\mH^3_{(X,Y,W)}(S_{(X,Y,W)})$ can be written uniquely as
$$
\sum_{i,j,k>0}\gfrac{\varphi_{ijk}}{W^i,X^j,(YZ)^k},
$$
where $\varphi_{ijk}\in\kappa(Z)$.
\end{cor}
\begin{proof}
$S_{(X,Y,W)}\simeq\kappa(Z)[X,YZ,W]_{(X,YZ,W)}$.
\end{proof}
\begin{cor}\label{cor:f}
Let $(f)$ be a non-zero prime ideal of $\kappa[Z,W]$ contained in 
$(Z,W)$ but not
containing $Z$ or $W$. An element $\Psi$ in $\mH^3_{(X,Y,f)}(S_{(X,Y,f)})$ can
be written as
$$
\Psi=\sum_{i,j>0}\gfrac{g_{ij}}{h_{ij},(XW)^i,(YZ)^j},
$$
where $g_{ij}\in\kappa[Z,W]$ and $0\neq h_{ij}\in\kappa[Z,W]$.  $\Psi=0$ if and
only if $g_{ij}\in h_{ij}\kappa[Z,W]_{(f)}$ for all $i,j$.
\end{cor}
\begin{proof}
Since $S_{(X,Y,f)}\simeq\kappa[Z,W]_{(f)}[XW,YZ]_{(XW,YZ)}$, there is an
isomorphism
$$
\mH^2_{(XW,YZ,f)}(S_{(X,Y,f)}\otimes_{\kappa[Z,W]_{(f)}}
\mH^1_{(f)}(\kappa[Z,W]_{(f)}))\simeq\mH^3_{(X,Y,f)}(S_{(X,Y,f)})
$$
\cite[(2.5)]{hu:pmzds} given by
$$
\sum_{i,j>0}\gfrac{1\otimes\gfrac{g_{ij}}{h_{ij}}}{(XW)^i,(YZ)^j}\mapsto
\sum_{i,j>0}\gfrac{g_{ij}}{h_{ij},(XW)^i,(YZ)^j}.
$$
The result follows from Proposition~\ref{prop:38555}. Moreover,
$\Psi=0$ if and only if $\gfrac{g_{ij}}{h_{ij}}=0$, equivalently
$g_{ij}\in h_{ij}\kappa[Z,W]_{(f)}$, for all $i,j$.
\end{proof}

The following lemma will be used in Section~\ref{sec:injresl}.
\begin{lem}\label{prop:onto}
Let $f\in\kappa[Z,W]$ be an irreducible polynomial in $(Z,W)$. For 
any $s,t>0$, there exist
$\ell>0$, $g\in\kappa[Z,W]$ and $h\in\kappa[Z,W]\setminus(f)$ such that
$$
\gfrac{g}{h,f^\ell}=\gfrac{1}{W^t,Z^s}
$$
in $\mH^2_{(Z,W)}(\kappa[Z,W]_{(Z,W)})$.
\end{lem}
\begin{proof}
We may assume $(f)\neq(W)$ to avoid the trivial case. Write
$$
f=f_0Z^u+f_1W^v
$$
for some $f_0\in\kappa[Z]\setminus(Z)$, 
$f_1\in\kappa[Z,W]\setminus(W)$ and $u,v>0$. Divide
$s$ by $u$:
\begin{eqnarray*}
s=uq+r & & \text{($0\leq q$ and $0\leq r< u$).}
\end{eqnarray*}
We choose $h$ to be $W^t$ and prove the lemma by
  induction
  on $\lceil t/v \rceil $, the smallest integer greater than or
  equal to $t/v$. In the case where $\lceil t/v\rceil =1$
  ({\em i.e.} $t\leq v$),
\[ f^{q+1}= (f_0Z^u + f_1W^v)^{q+1} = f_0^{q+1} Z^{u(q+1)} +w W^t, \]
for some $w \in \kappa[Z,W]$. The following can be computed
using Example~\ref{example2.2}:
$$
\gfrac{f_0^{q+1}Z^{u-r}}{W^t,f^{q+1}}
=\gfrac{f_0^{q+1}Z^{u-r}}{W^t,f_0^{q+1}Z^{uq+u}}=\gfrac{1}{W^t,Z^s}.
$$
Assume the lemma holds for $\lceil t/v\rceil=n$. For the case $\lceil t/v\rceil= n+1$,
let
\begin{eqnarray*}
F&=&\sum_{i=0}^{n-1}(f_0Z^u)^i(-f_1W^v)^{n-i-1},\\
G&=&\sum_{j=0}^q\binom{q+1}{j}(f_0Z^u)^{nj}(-f_1W^vF)^{q-j}.
\end{eqnarray*}
Then
\begin{eqnarray*}
f((f_0Z^u)^n-f_1W^vF)&=&(f_0Z^u)^{n+1}-(-f_1W^v)^{n+1},\\
((f_0Z^u)^n-f_1W^vF)^{q+1}&=&(f_0Z^u)^{n(q+1)}-f_1W^vFG,
\end{eqnarray*}
and
\begin{eqnarray*}
\gfrac{f_0^{q+1}Z^{u-r}}{W^t,f^{q+1}}
&=&\gfrac{f_0^{q+1}Z^{u-r}((f_0Z^u)^n-f_1W^vF)^{q+1}}{W^t,((f_0Z^u)^{n+1}-(-f_1W^v)^{n+1})^{q+1}}\\
&=&\gfrac{f_0^{q+1}Z^{u-r}((f_0Z^u)^n-f_1W^vF)^{q+1}}{W^t,(f_0Z^u)^{(n+1)(q+1)}}\\
&=&\gfrac{f_0^{q+1}Z^{u-r}(f_0Z^u)^{n(q+1)}}{W^t,(f_0Z^u)^{(n+1)(q+1)}}-
\gfrac{f_0^{q+1}Z^{u-r}f_1W^vFG}{W^t,(f_0Z^u)^{(n+1)(q+1)}}\\
&=&\gfrac{1}{W^t,Z^s}-
\gfrac{f_0^{-n(q+1)}Z^{u-r}f_1FG}{W^{t-v},Z^{u(n+1)(q+1)}}.
\end{eqnarray*}
Since $\lceil(t-v)/v\rceil=n$, there exist $\ell_0>0$ and 
$g_0\in\kappa[Z,W]$ such that
$$
\gfrac{g_0}{W^{t-v},f^{\ell_0}}=\gfrac{1}{W^{t-v},Z^{u(n+1)(q+1)}}.
$$
We get the required elements:
\begin{eqnarray*}
&
&\gfrac{f^{\ell_0}f_0^{q+1}Z^{u-r}+f^{q+1}g_0f_0^{-n(q+1)}Z^{u-r}W^vf_1FG}{W^t,f^{\ell_0+q+1}}\\
&=&
\gfrac{f_0^{q+1}Z^{u-r}}{W^t,f^{q+1}}+\gfrac{g_0f_0^{-n(q+1)}Z^{u-r}f_1FG}{W^{t-v},f^{\ell_0}}
\\
&=&
\gfrac{f_0^{q+1}Z^{u-r}}{W^t,f^{q+1}}+\gfrac{f_0^{-n(q+1)}Z^{u-r}f_1FG}{W^{t-v},Z^{u(n+1)(q+1)}}
=\gfrac{1}{W^t,Z^s}.
\end{eqnarray*}

\end{proof}


\section{injective hulls}\label{sec:injmod}


In this section, we study the module structure of an injective hull $E(\q)$ of
$A/\fQ$ for each prime ideal $\fQ$ of $A$ generated by $X,Y$ and a
prime ideal $\q$ of $\kappa[Z,W]$ contained in $(Z,W)$. We use the 
following two
well-known constructions for injective hulls.
\begin{lem}\label{lem:1}
Let $R$ be a Noetherian ring, $\fQ\subset\fP$ be prime ideals of $R$ 
and $E(R/\fQ)$ be an injective
hull of $R/\fQ$. Then $E(R/\fQ)$ is an $R_\fP$-module and it is an 
injective hull of
$(R/\fQ)_\fP$ over $R_\fP$.
\end{lem}
\begin{lem}\label{lem:2}
Let $R$ be a Noetherian ring, $I$ be an ideal of $R$, $\fP$ be a 
prime ideal of $R$ containing $I$ and
$E(R/\fP)$ be an injective hull of $R/\fP$. Then, as an $R/I$-module, 
$\Hom_R(R/I,E(R/\fP))$ is an
injective hull of $R/\fP$.
\end{lem}

Let $\q=(f_1,\cdots,f_m)$ be a prime ideal of $\kappa[Z,W]$ contained 
in $(Z,W)$
and $\fQ$ be the prime ideal of $S$ generated by 
$X,Y,f_1,\cdots,f_m$. We denote
by $\fQ$ also the element of $\Spec A$, $\Spec S_\fQ$ and $\Spec 
S_{(X,Y,Z,W)}$,
which canonically embed into $\Spec S$. Recall that
$\mH^{\height\fQ}_{\fQ}(S_{\fQ})$ is an injective hull of $S/\fQ$, as $S$ is
a Gorenstein ring.
\begin{defn}
$$
E(f_1,\cdots,f_m):=E(\q):=\{\omega\in\mH^{\height\fQ}_{\fQ}(S_{\fQ})|
XW\omega=YZ\omega\}
$$
\end{defn}
By Lemma~\ref{lem:1}, $\mH^{\height\fQ}_{\fQ}(S_{\fQ})$ as an
$S_{(X,Y,Z,W)}$-module is also an injective hull of
$S_{(X,Y,Z,W)}/\fQ$. By Lemma~\ref{lem:2}, with the $A$-module
structure via the bijection
$$
E(\q)\simeq\Hom_{S_{(X,Y,Z,W)}}(A,\mH^{\height\fQ}_{\fQ}(S_{\fQ})),
$$
$E(\q)$ is an injective hull of $A/\fQ$. Next, we describe elements 
in $E(\q)$ using certain maps
$\Omega_{\q}^n$. If $\q$ is principal,
$$
\Omega^n_\q\colon\kappa(Z,W)\to\mH^{\height\fQ}_{\fQ}(S_{\fQ})
$$
is a $\kappa[Z,W]_\q$-linear map. If
$\q=(Z,W)$,
$$
\Omega^n_\q\colon\mH^2_{(Z,W)}(\kappa[Z,W]_{(Z,W)})\to
\mH^4_{(X,Y,Z,W)}(S_{(X,Y,Z,W)})
$$
is a $\kappa$-linear map.
$\Omega^n_\q$ is defined to be zero for $n<0$ and is defined below 
for $n\geq 0$.

\begin{defn}\label{def:1522}
Let $n\geq 0$ and $\q=(f)$ be a prime ideal of $\kappa[Z,W]$ 
contained in $(Z,W)$. Given
$s\in\mathbb Z$, $g\in\kappa[Z,W]$ and $h\in\kappa[Z,W]\setminus\q$, we define
\begin{eqnarray*}
\Omega^n_\q(\frac{g}{h}):=
\dfrac{g}{h}\displaystyle\sum_{i=0}^{n}\gfrac{1}{(XW)^{i+1},(YZ)^{n+1-i}},
& & \text{if $(f)=(0)$;}
\end{eqnarray*}
and
$$
\Omega^n_\q(\frac{g}{hf^s}):=
\begin{cases}
\dfrac{g}{h}\displaystyle\sum_{i=0}^{n}\gfrac{1}{f^sZ^{n+1-i},(XW)^{i+1},Y^{n+1-i}},
& \text{if $(f)=(Z)$;}\\
\dfrac{g}{h}\displaystyle\sum_{i=0}^{n}\gfrac{1}{f^sW^{i+1},X^{i+1},(YZ)^{n+1-i}},
& \text{if $(f)=(W)$;}\\
\dfrac{g}{h}\displaystyle\sum_{i=0}^{n}\gfrac{1}{f^s,(XW)^{i+1},(YZ)^{n+1-i}},
& \text{if $(f)\neq(0)$, $(Z)$ or $(W)$.}
\end{cases}
$$
The $\kappa$-linear map $\Omega^n_{(Z,W)}$ is defined by
$$
\Omega^n_{(Z,W)}\gfrac{1}{Z^u,W^v}:=
\sum_{i=0}^{n}\gfrac{1}{Z^{n+1-i+u},W^{i+1+v},X^{i+1},Y^{n+1-i}},
$$
where $u,v>0$.
\end{defn}
If $\q$ is principal, $\Omega^n_\q$ is independent of the choice of a 
generator $f$.
We use also the notation $\Omega^n_f:=\Omega^n_{\q}$. If $\q=(Z,W)$, 
we use also the
notation $\Omega^n_{Z,W}:=\Omega^n_\q$.
The following facts are not hard to check. Details are left to the reader.
\begin{prop}\label{prop:888311}
Let $n\geq 0$.
\begin{enumerate}
\item\label{it:12394} For $\varphi\in\kappa(Z,W)$, 
$\Omega_0^n(\varphi)\neq 0$ if and only if
           $\varphi\neq 0$.
\item $\Omega^n_Z(Z^s)\neq 0$ (resp. $\Omega^n_W(W^s)\neq 0$) if and only if
           $s\leq n$. In the $\kappa$-vector space $E(Z)$ (resp. 
$E(W)$), elements of the form
           $\Omega^n_Z(Z^sW^t)$ (resp. $\Omega^n_W(Z^tW^s)$), where
           $s\leq n$ and $t\in\mathbb Z$, are linearly independent.
\item For non-zero $f$ with $Z,W\not\in(f)$,
           $\Omega^n_f(f^s)\neq 0$ if and only if $s<0$.
\end{enumerate}
\end{prop}

For any $\varphi\in\kappa(Z,W)$,
$$
XW\Omega^n_f\varphi=YZ\Omega^n_f\varphi=\Omega^{n-1}_f\varphi.
$$
Therefore $\Omega^n_f$ has image in $E(f)$. The multiplication
by $XW$ (equals $YZ$ in $A$) takes elements of $E_n(\varphi)$ into
$E_{n-1}(\varphi)$.
\begin{example}\label{ex:E(f)}
For any non-zero $\varphi$ in $\kappa(Z,W)$,
$$
(XW)^n \Omega_0^n(\varphi) = (YZ)^n \Omega_0^n(\varphi)=  \Omega_0^0(\varphi)
$$
and it is non-zero in $E(f)$ by Proposition~\ref{prop:888311}~(\ref{it:12394}).
\end{example}

Next, we explain the structure of $E(Z,W)$.
For any  $\varphi\in\mH^2_{(Z,W)}(\kappa[Z,W]_{(Z,W)})$,
$$
XW\Omega^n_{Z,W}\varphi=YZ\Omega^n_{Z,W}\varphi=\Omega^{n-1}_{Z,W}\varphi.
$$
Therefore $\Omega^n_{Z,W}$ has image in $E(Z,W)$.  Note that 
$\Omega^n_{Z,W}$ is
not $\kappa[Z,W]$-linear for $n\geq 0$. For instance, 
$\gfrac{1}{Z,W}$ is annihilated by
$W^{n+1}$ but
$$
W^{n+1}\Omega^n_{Z,W}\gfrac{1}{Z,W}=\gfrac{1}{Z^2,W,X^{n+1},Y}\neq 0.
$$

\begin{defn}
For $s,t\in\mathbb Z$, we choose $u,v>0$
with $u+s,v+t\geq 0$ and define the notation
$$
\Omega^n(Z^sW^t):=Z^{u+s}W^{v+t}\Omega^n_{Z,W}\gfrac{1}{Z^u,W^v}.
$$
\end{defn}
This definition is independent of the choice of $u$ and $v$, indeed,
\begin{equation}\label{eq:341112}
\Omega^n(Z^sW^t)=\sum_{i=0}^{n}\gfrac{1}{Z^{n+1-i-s},W^{i+1-t},X^{i+1},Y^{n+1-i}}.
\end{equation}
In general, $\Omega^n(Z^sW^t)$ does not equal to
$\Omega^n_{Z,W}\gfrac{1}{Z^{-s},W^{-t}}$. For instance,
$$
\Omega^1_{Z,W}\gfrac{1}{Z,W^{-1}}=0,
$$
but
$$
\Omega^1(Z^{-1}W)=W^2\Omega^1_{Z,W}\gfrac{1}{Z,W}\neq 0.
$$

\begin{prop}\label{prop:058844}
$\Omega^n(Z^sW^t)\neq 0$ if and only if $n\geq\max\{0,s,t,s+t\}$. The 
non-trivial
$\Omega^n(Z^sW^t)$ form a basis for the $\kappa$-vector space $E(Z,W)$.
\end{prop}
\begin{proof}
Let $n\geq 0$. It is clear from the definition that 
$\Omega^n(Z^sW^t)=0$ if one of
$s$, $t$, and $t+s$ is greater than $n$.
We show first that the elements of the form $\Omega^n(Z^sW^t)$ 
generate $E(Z,W)$. Let
$$
\Psi=\sum_{i,j,k,\ell\geq 1}\gfrac{a_{ijk\ell}}{Z^i,W^j,X^k,Y^\ell}\in
\mH^4_{(X,Y,Z,W)}(S_{(X,Y,Z,W)}),
$$
where $a_{ijk\ell}\in\kappa$. Assume that $\Psi\in E(Z,W)$, that is,
$XW\Psi=YZ\Psi$ or
$$
\sum_{i,j,k,\ell\geq 1}\gfrac{a_{ijk\ell}}{Z^i,W^{j-1},X^{k-1},Y^\ell}=
\sum_{i,j,k,\ell\geq 1}\gfrac{a_{ijk\ell}}{Z^{i-1},W^j,X^k,Y^{\ell-1}}.
$$
Comparing coefficients, we get
$$
a_{i(j+1)(k+1)\ell}=a_{(i+1)jk(\ell+1)}
$$
for $i,j,k,\ell\geq 1$. For $i',j',k',$ or $\ell'$ less than $1$, if 
there exist
$i,j,k,\ell\geq 1$ such that $i+j=i'+j'$, $k+\ell=k'+\ell'$ and $j-k=j'-k'$, we
define
$$
a_{i'j'k'\ell'}:=a_{ijk\ell};
$$
otherwise we define $a_{i'j'k'\ell'}:=0$. Then
\begin{eqnarray*}
\Psi&=&\sum_{\stackrel{\scriptstyle{\ell,m\in\mathbb Z}}{n\geq 0}}
a_{\ell m1(n+1)}\sum_{i=0}^n\gfrac{1}{Z^{\ell-i},W^{m+i},X^{1+i},Y^{n+1-i}}\\
&=&\sum_{\stackrel{\scriptstyle{\ell,m\in\mathbb Z}}{n\geq 0}}
a_{\ell m1(n+1)}\Omega^n(Z^{n+1-\ell}W^{1-m}).
\end{eqnarray*}

Now we show that those $\Omega^n(Z^sW^t)$ with $n\geq max\{0,s,t,s+t\}$ are
linearly independent over $\kappa$. We study a linear combination of
$\Omega^n(Z^sW^t)$:
\begin{equation}\label{lincomb}
\sum_{n\geq max\{0,s,t,s+t\}} a_{nst}\Omega^n(Z^sW^t),
\end{equation}
where $a_{nst}\in\kappa $. Setting $i$ equal to $n$ and $n-s$ 
respectively for the
expression (\ref{eq:341112}) of $\Omega^n(Z^sW^t)$, we have
\begin{eqnarray*}
\lefteqn{\gfrac{1}{Z^{n+1-i-s},W^{i+1-t},X^{i+1},Y^{n+1-i}}=} \\
&&
\begin{cases}
  \gfrac{1}{Z^{1-s},W^{n-t+1},X^{n+1},Y}, & \text{ if $i=n$; } \\ \\
  \gfrac{1}{Z,W^{n-s-t+1},X^{n-s+1},Y^{s+1}}, & \text{ if $i=n-s$ }.
\end{cases}\end{eqnarray*}
We note that
\begin{eqnarray*}
\gfrac{1}{Z^{1-s},W^{n-t+1},X^{n+1},Y} = 0, && \text{if $s > 0$}
\end{eqnarray*}
and
\begin{eqnarray*}
\gfrac{1}{Z,W^{n-s-t+1},X^{n-s+1},Y^{s+1}} =0 , && \text{if $s<0$.}
\end{eqnarray*}

For given $n, s, t$ in the summation in (\ref{lincomb}),
$\gfrac{1}{Z^{1-s},W^{n-t+1},X^{n+1},Y}$ occurs if $s <0$.
By Corollary~\ref{cor:ZW}, $\gfrac{1}{Z^i,W^j,X^k,Y^{\ell}}$ are linearly
independent for all $i,j,k,\ell >0$ in
$H^4_{(X,Y,Z,W)}(S_{(X,Y,Z,W)})$. This implies that if there exist
$n',s',t'$ such that
\[ \gfrac{1}{Z^{n'+1-i'-s'},W^{i'+1-t'},X^{i'+1},Y^{n'+1-i'}} =
\gfrac{1}{Z^{1-s},W^{n-t+1},X^{n+1},Y}\] for some $i' \in \{0, \dots,
n' \}$, then $i'=n'=n$, $s'=s$ and $t'=t$. Thus, for fixed $n,s,t$ 
with $s < 0$,
the generalized fraction $\gfrac{1}{Z^{1-s},W^{n-t+1},X^{n+1},Y}$ occurs in
(\ref{lincomb}) exactly once with the coefficient $a_{nst}$. Similarly for
$\gfrac{1}{Z,W^{n-s-t+1},X^{n-s+1},Y^{s+1}}$ with $s \geq 0$.
Therefore, if
\[ \sum_{n \geq max\{0,s,t,s+t\}} a_{nst} \Omega^n(Z^sW^t) =0, \]
then $a_{nst}=0$ for all $n,s,t$ by
Corollary~\ref{cor:ZW} again. Hence, $\Omega^n(Z^s W^t)$ are linearly 
independent.
\end{proof}

The $A$-module structure of $E(Z,W)$ is clear:
For $s_1,t_1\geq 0$ and $s_2,t_2\in\mathbb Z$, we have
\begin{eqnarray}
Z^{s_1}W^{t_1}\Omega^n(Z^{s_2}W^{t_2})&=&\Omega^n(Z^{s_1+s_2}W^{t_1+t_2}),\label{eq:87761}\\
X^{t_1}Y^{s_1}\Omega^n(Z^{s_2}W^{t_2})&=&\Omega^{n-t_1-s_1}(Z^{s_2-s_1}W^{t_2-t_1}).\label{eq:87762}
\end{eqnarray}
For arbitrary $\varphi\in A$, $n\geq 0$ and $s_2+t_2\leq n$, we choose
$f\in\kappa[X,Y,Z,W]$ such that
$\varphi-f\in(X,Y)^{n+1}+(Z,W)^{n-s_2-t_2+1}$, then
$$
\varphi\Omega^n(Z^{s_2}W^{t_2})=f\Omega^n(Z^{s_2}W^{t_2}).
$$
Replaced $\varphi$ by $f$, we can use the equalities (\ref{eq:87761}),
(\ref{eq:87762}) and $\kappa$-linearity to multiply 
$\Omega^n(Z^{s_2}W^{t_2})$ by
$\varphi$.

The multiplication by $X$ ({\em resp.} $Y$) takes $\Omega^n(Z^sW^t)$
to $\Omega^{n-1}(Z^sW^{t-1})$ ({\em resp.}
$\Omega^{n-1}(Z^{s-1}W^t)$) in $E_{n-1}(Z,W)$. Alternative to 
Example~\ref{ex:E(f)},
multiplications by variables can take non-zero elements in $E(Z,W)$
to non-zero elements in $E_{0}(Z,W)$.
\begin{example}\label{ex:E(Z,W)}
By Proposition
\ref{prop:058844}, for any non-zero $\Omega^n(Z^sW^t)$, there exist
$n_1$ and $n_2$ with $n_1+n_2 =n$ such that
$$
X^{n_1}Y^{n_2}\Omega^n(Z^sW^t)= \Omega^0(Z^{s-n_2}W^{t-n_1}) \neq 0.
$$
Note that $\Omega^0(Z^sW^t)$ obtained by multiplying
$\Omega^n(Z^sW^t)$ by $(XW)^n$ may be zero if
$s$ or $t$ is positive.
\end{example}

The computations in Examples \ref{ex:E(f)} and
\ref{ex:E(Z,W)} will be used in proving our main result Theorem \ref{thm:main}.

Divisions by $X$, $Y$, $Z$ and $W$ can be defined as well. In general,
for $i,j,k,l\in\mathbb Z$, let $X^iY^jZ^kW^l$ be the
$\kappa$-linear operator on $E(Z,W)$ satisfying
$$
X^iY^jZ^kW^l\Omega^n(Z^sW^t)=\Omega^{n-i-j}(Z^{s+k-j}W^{t+l-i}).
$$
Using the above description of the $A$-module structure of $E(Z,W)$, 
One can check that
this operator is $A$-linear.
\begin{defn}\label{defn:43890}
Let $\q=(f_1,\cdots,f_m)$ be a prime ideal of $\kappa[Z,W]$ contained 
in $(Z,W)$. We
define $E_n(\q)$, denoted also by $E_n(f_1,\cdots,f_m)$, to be the
$\kappa[Z,W]_{(Z,W)}$-submodule of $E(\q)$ generated by the image of 
$\Omega^n_\q$.
\end{defn}
$E_n(f)$ consists of elements of the form $\Omega^n_f\varphi$. 
$E_n(Z,W)$ consists of
elements of the form $Z^sW^t\Omega^n_{Z,W}\varphi$. Note that powers 
of $Z$ and $W$ are
necessary to represent elements of $E_n(Z,W)$. For instance,
$$
\gfrac{1}{Z,W,X^2,Y}=ZW^2\Omega^1_{Z,W}\gfrac{1}{Z,W}
$$
does not equal to $\Omega^n_{Z,W}\varphi$ for any 
$\varphi\in\mH^2_{(Z,W)}(\kappa[Z,W]_{(Z,W)})
$.
\begin{prop}\label{prop:34312}
$E(\q)=\oplus_n E_n(\q)$ as $\kappa[Z,W]_{(Z,W)}$-modules.
\end{prop}
\begin{proof}
We prove the proposition in five cases.
\begin{trivlist}
\item{Case 1. $\q=(Z,W)$.} Already shown in Proposition~\ref{prop:058844}.
\item{Case 2. $\q=(0)$.} Let
$$
\Psi=\sum_{i,j\geq 
1}\gfrac{\varphi_{ij}}{(XW)^i,(YZ)^j}\in\mH^2_{(X,Y)}(S_{(X,Y)}),
$$
where $\varphi_{ij}\in\kappa(Z,W)$, be an element of $E(0)$. From the identity
$XW\Psi=YZ\Psi$, we get
$$
\varphi_{i(j+1)}=\varphi_{(i+1)j}
$$
for $i,j\geq 1$. We have the expression
$$
\Psi=\sum_{n\geq 
0}\left(\sum_{i+j=n+2}\gfrac{\varphi_{ij}}{(XW)^i,(YZ)^j}\right)=
\sum_{n\geq 0}\Omega^n_0(\varphi_{1(n+1)}),
$$
since for any $i,j$ with $i+j=n+2$, it is clear that
$\varphi_{ij}=\varphi_{i(n+2-i)}=\varphi_{1(n+1)}$. This shows
$E(0)=\sum_n E_n(0)$.

An element
$$
\Psi_n=\sum\gfrac{\varphi_{ij}}{(XW)^i,(YZ)^j}\in E_n(0),
$$
where $\varphi_{ij}\in\kappa(Z,W)$, satisfies $i+j=n+2$. If
$\sum\Psi_n=0$, by Corollary~\ref{cor:0}, $\Psi_n=0$ for all $n$. This shows
$E(0)=\oplus_n E_n(0)$.
\item{Case 3. $\q=(Z)$.} Let
$$
\Psi=\sum_{i,j,k\geq 1}\gfrac{\varphi_{ijk}}{Z^i,(XW)^j,Y^k}
\in\mH^3_{(X,Y,Z)}(S_{(X,Y,Z)})
$$
where $\varphi_{ijk}\in\kappa(W)$, be an element of $E(Z)$. From the identity
$XW\Psi=YZ\Psi$, we get
$$
\varphi_{i(j+1)k}=\varphi_{(i+1)j(k+1)}
$$
for $i,j,k\geq 1$. For $i',j'$, or $k'$ less than $1$, if there exist
$i,j,k\geq 1$ such that $i+j=i'+j'$ and $j+k=j'+k'$, we define
$$
\varphi_{i'j'k'}:=\varphi_{ijk};
$$
otherwise we define $\varphi_{i'j'k'}:=0$. Then
$$
\Psi=\sum_{\stackrel{\scriptstyle{m\in\mathbb Z}}{n\geq 0}}\varphi_{m1(n+1)}
\sum_{i=0}^n
\gfrac{1}{Z^{m-i},(XW)^{i+1},Y^{n+1-i}}=
\sum_{\stackrel{\scriptstyle{m\in\mathbb Z}}{n\geq 0}}\Omega^n_Z
(\frac{\varphi_{m1(n+1)}}{Z^{m-n-1}}).
$$
This shows $E(Z)=\sum_n E_n(Z)$.

An element
$$
\Psi_n=\sum_{i=0}^{n}\gfrac{\varphi_{ijk}}{Z^i,(XW)^j,Y^k}\in E_n(Z),
$$
where $\varphi_{ijk}\in\kappa(W)$, satisfies $j+k=n+2$. If
$\sum\Psi_n=0$, by Corollary~\ref{cor:Z}, $\Psi_n=0$ for all $n$. This shows
$E(Z)=\oplus_n E_n(Z)$.
\item{Case 4. $\q=(W)$.} Similar to Case 3.
\item{Case 5. $\q=(f)$ not equal to $(Z)$ or $(W)$.} Let
$$
\Psi=\sum_{i,j\geq 1}\gfrac{g_{ij}}{h_{ij},(XW)^i,(YZ)^j}
\in\mH^3_{(X,Y,f)}(S_{(X,Y,f)}),
$$
where $g_{ij}\in\kappa[Z,W]$ and $0\neq h_{ij}\in\kappa[Z,W]$, be an element in
$E(f)$. Multiplying the numerators and denominators by $h_{ij}$'s, we 
may assume
that all $h_{ij}$ equals a fixed $h\in\kappa[Z,W]$. From the identity
$XW\Psi=YZ\Psi$, we get
$$
g_{i(j+1)}-g_{(i+1)j}\in h\kappa[Z,W]_{(f)}
$$
for $i,j\geq 1$. Hence
$$
\Psi=\sum_{n\geq 0}\left(\sum_{i+j=n+2}
\gfrac{g_{1(n+1)}}{h,(XW)^i,(YZ)^j}\right)=
\sum_{n\geq 0}\Omega^n_f(\frac{g_{1(n+1)}}{h}).
$$
This shows $E(f)=\sum_n E_n(f)$.

An element
$$
\Psi_n=\sum_{i,j\geq 1}\gfrac{g_{ij}}{h,(XW)^i,(YZ)^j}\in E_n(f),
$$
where $g_{ij}\in\kappa[Z,W]$ and $0\neq h\in\kappa[Z,W]$, satisfies
$i+j=n+2$. If $\sum\Psi_n=0$, by Corollary~\ref{cor:f}, $\Psi_n=0$ 
for all $n$. This
shows $E(f)=\oplus_n E_n(f)$.
\end{trivlist}
\end{proof}
$E_n(\q)$ is not a $\kappa[X,Y,Z,W]$-module. In fact, $XWE_n(\q)=E_{n-1}(\q)$.
\begin{prop}\label{prop:086661}
Let $\q$ be a prime ideal of $\kappa[Z,W]$ contained in $(Z,W)$. Then
$E_0(\q)=0\colon_{E(\q)}(X,Y)$.
\end{prop}
It is clear that $E_0(\q)$ is annihilated by $(X,Y)$ for all $\q$. To 
prove the proposition,
it remains to show that $E_0(\q)$ contains all the elements 
annihilated by $(X,Y)$.
For a prime ideal $(f)$ of $\kappa[Z,W]$ contained in $(Z,W)$, we denote
$$
f^\triangle:=
\begin{cases}
X,&\text{if $(f)=(W)$};\\
Y,&\text{if $(f)=(Z)$};\\
XW,&\text{otherwise}.
\end{cases}
$$
The multiplication on $E(f)$ by $Z$ (resp. $W$) is an isomorphism if 
$Z\not\in(f)$
(resp. $W\not\in(f)$). Since elements of $E(f)$ are annihilated by 
$XW-YZ$, an element
of $E(f)$ is annihilated by $f^\triangle$ if and only if it is 
annihilated by $X$ and
$Y$. For instance, let $f=Z+W$ and $\Psi\in E(f)$. Then $f^\triangle=XW$ and
$(XW-YZ)\Psi=0$. If $X\Psi=Y\Psi=0$, then clearly 
$f^\triangle\Psi=0$. Conversely if
$f^\triangle\Psi=0$, then $X\Psi=0$, since the multiplication by $W$ 
is an isomorphism.
$YZ\Psi=0$ as well, since $YZ\Psi=f^\triangle\Psi$. Now $Y\Psi=0$, because the
multiplication by $Z$ is an isomorphism.

For Proposition~\ref{prop:086661}, what we need to prove is the following.
\begin{prop}\label{prop:876111}
$E_0(f)$ contains all the elements of $E(f)$ annihilated by 
$f^\triangle$. $E_0(Z,W)$
contains all the elements of $E(Z,W)$ annihilated by $(X,Y)$.
\end{prop}
\begin{proof}
Let
$$
\Psi=\sum\Omega^n_0(\varphi_n)
$$
be an element of $E(0)$ annihilated by $XW$, where 
$\varphi_n\in\kappa(Z,W)$. Then
$$
\sum\Omega^{n-1}_0(\varphi_n)=XW\Psi=0.
$$
By Proposition \ref{prop:34312}, $\Omega^{n-1}_0(\varphi_n)=0$ for 
all $n$, which implies $\varphi_n=0$ for $n\geq 1$ by Proposition~\ref{prop:888311}. 
Therefore
$$
\Psi=\Omega^0_0(\varphi_0)\in E_0(0).
$$

Now assume that $(f)\neq 0$, $(Z)$ or $(W)$. Let
$$
\Psi=\sum\Omega^n_f(\frac{g_n}{h_n}f^{s_n})
$$
be an element of $E(f)$ annihilated by $XW$, where
$g_n,h_n\in\kappa[Z,W]\setminus(f)$.  Then
$$
\sum\Omega^{n-1}_f(\frac{g_n}{h_n}f^{s_n})=XW\Psi=0.
$$
By Propositions~\ref{prop:34312} and \ref{prop:888311},
$\Omega^{n-1}_f(\frac{g_n}{h_n}f^{s_n})=0$ for all
$n$ and this implies $s_{n}\geq 0$ for $n\geq 1$. Therefore
$\Omega^n_f(\frac{g_n}{h_n}f^{s_n})=0$
for all $n\geq 1$ and
$$
\Psi=\Omega^0_f(\frac{g_0}{h_0}f^{s_0})\in E_0(f).
$$

Let
$$
\Psi=\sum\Omega^n_Z(\frac{g_n}{h_n}Z^{s_n})
$$
be an element of $E(Z)$ annihilated by $Y$, where
$g_n,h_n\in\kappa[Z,W]\setminus(Z)$. Then
$$
\sum_{n\geq 1}\Omega^{n-1}_Z(\frac{g_n}{h_n}Z^{s_{n}-1})=Y\Psi=0.
$$
Using Propositions~\ref{prop:34312} and \ref{prop:888311} again, a similar
argument as in the previous cases shows
$\Omega^{n-1}_Z(\frac{g_n}{h_n}Z^{s_{n}-1})=0$ for all $n$ and
$s_{n}-1>n-1$ for $n\geq 1$. Therefore
$\Omega^n_Z(\frac{g_n}{h_n}Z^{s_n})=0$
for all $n\geq 1$ and
$$
\Psi=\Omega^0_Z(\frac{g_0}{h_0}Z^{s_0})\in E_0(Z).
$$

Similarly, we see that $E_0(W)$ contains all the elements of $E(W)$ 
annihilated by $X$.

Let
$$
\Psi=\sum_{n\geq\max\{0,s,t,s+t\}}a_{nst}\Omega^n(Z^sW^t)
$$
be an element of $E(Z,W)$ annihilated by $X$ and $Y$, where $a_{nst}\in\kappa$.
Then
$$
\sum_{n\geq\max\{1,s+1,t,s+t\}}a_{nst}\Omega^{n-1}(Z^sW^{t-1})=X\Psi=0.
$$
By Proposition~\ref{prop:058844}, the coefficient $a_{nst}=0$, if 
$n\geq\max\{1,s+1,t,s+t\}$. For
$n\geq 1$, possible non-trivial coefficients are those $a_{nnt}$ with 
$t\leq 0$. Similarly,
$Y\Psi=0$ implies that possible non-trivial coefficients are those 
$a_{nsn}$ with $s\leq 0$.
Therefore $a_{nst}=0$ for $n\geq 1$ and
$$
\Psi=\sum_{0\geq\max\{s,t,s+t\}}a_{0st}\Omega^0(Z^sW^t)\in E_0(Z,W).
$$\end{proof}


\section{an injective resolution}
\label{sec:injresl}


In this section, we construct explicitly an injective resolution of 
$A/\p$ using the injective
modules given in Section~\ref{sec:injmod}. The coboundary maps of the 
injective resolution
involve multiplications and divisions by elements of $A$ and  certain 
maps $d^0_f$,
$d^1_f$ appeared in a residual complex.

A prime ideal $(f)$ is also generated by $gf$ for any invertible 
element $g$. For convenience, we
use the notation $\oplus_{f\neq 0}E(f)$ for the direct sum of modules 
$E(f)$ indexed by
the ideals generated by the irreducible polynomial $f\in\m=(Z,W)$; that is,
$$
\underset{\scriptscriptstyle f\neq 0}{\oplus}E(f):=
\underset{\scriptscriptstyle \q\neq(0),\m}{\oplus}E(\q),
$$
where $f$ ranges over irreducible polynomials contained in $(Z,W)$. 
We use the notation
$\sum_f$ for representing elements in $\oplus_{f\neq 0}E(f)$.

Recalling the notation $f^\triangle$ defined in 
Section~\ref{sec:injmod}, we have the
following exact sequence by Propositions~\ref{prop:086661} and 
\ref{prop:876111}.
\begin{equation}\label{eq:388875}
0\to\underset{\scriptscriptstyle f\neq 0}{\oplus}E_0(f)
\to
\underset{\scriptscriptstyle f\neq 0}{\oplus}E(f)
\xrightarrow{\oplus f^{\triangle}}
\underset{\scriptscriptstyle f\neq 0}{\oplus}E(f)
\to 0.
\end{equation}

Now we define $d^0_f$ using Corollary \ref{cor:0}.
\begin{defn}
For an irreducible polynomial $f\in(Z,W)$, we define
$$
d^0_f\colon\mH^2_{(X,Y)}(S_{(X,Y)})\to\mH^3_{(X,Y,f)}(S_{(X,Y,f)})
$$
to be the  map
$$
\sum_{i,j>0}\gfrac{g_{ij}/h_{ij}}{(XW)^i,(YZ)^j}\mapsto
\sum_{i,j>0}\gfrac{g_{ij}}{h_{ij}W^iZ^j,X^i,Y^j},
$$
where $g_{ij}\in\kappa[Z,W]$ and $0\neq h_{ij}\in\kappa[Z,W]$.
\end{defn}
We note that $d^0_f$ induces a restriction (by abusing the notation)
$$
d^0_f\colon E(0)\to E(f)
$$
since $d^0_f\left(\Omega^n_0(g/h)\right)=\Omega^n_f(g/h)$
and $d^0_f(E_n(0))\subset E_n(f)$. The product
$$
\prod_{f\neq 0}d^0_f\colon E(0)\to\prod_{f\neq 0}E(f)
$$
has image in $\oplus_{f\neq 0}E(f)$.
\begin{defn}
We define
$$
d^0\colon E(0)\to\underset{\scriptscriptstyle f\neq 0}{\oplus}E(f),
$$
where $f$ ranges over irreducible polynomials in $(Z,W)\kappa[Z,W]$, to be the
$A$-linear map
$$
\Omega^n_0(\frac{g}{h})\mapsto\sum_{f\neq 0} \Omega^n_f(\frac{g}{h}),
$$
where $g \in \kappa[Z,W]$ and $0\neq h\in\kappa[Z,W]$.
\end{defn}

Now we define $d^1_f$ using Corollaries \ref{cor:Z}, \ref{cor:W} and 
\ref{cor:f}.
\begin{defn}
For an irreducible polynomial $f\in(Z,W)$, we define
$$
d^1_f\colon
\mH^3_{(X,Y,f)}(S_{(X,Y,f)})\to\mH^4_{(X,Y,Z,W)}(S_{(X,Y,Z,W)})
$$
to be the map
$$
\sum_{i,j,k>0}\gfrac{g_{ijk}/h_{ijk}}{f^i,X^j,Y^k}\mapsto
\sum_{i,j,k>0}\gfrac{g_{ijk}}{h_{ijk},f^i,X^j,Y^k},
$$
where $g_{ijk}$ and $h_{ijk}$ are in $\kappa[Z,W]$ and $h_{ijk}$ has 
no factor $f$.
\end{defn}
For instance,
\begin{eqnarray}
d^1_W\Omega^n_W(Z^sW^t)&=&\Omega^n(Z^sW^t)\label{eq:88068}\\
d^1_Z\Omega^n_Z(Z^sW^t)&=&-\Omega^n(Z^sW^t).\nonumber
\end{eqnarray}

\begin{lem}
$d^1_f(E_n(f))\subset E_n(Z,W)$
\end{lem}
\begin{proof}
Let $g,h\in\kappa[Z,W]\setminus(f)$.
\begin{trivlist}
\item[Case 1. $(f)\neq(Z)$ or $(W)$.] For $j>0$, the elements $hW^{n+1}Z^{n+1}$
and $f^j$ form a system of parameters for $k[Z,W]_{(Z,W)}$ , so there exist
$\alpha_{ij}\in\kappa[Z,W]_{(Z,W)}$ and $s,t>n+1$ such that
$$
\begin{cases}
Z^s=\alpha_{11}hW^{n+1}Z^{n+1}+\alpha_{12}f^j,\\
W^t=\alpha_{21}hW^{n+1}Z^{n+1}+\alpha_{22}f^j.
\end{cases}
$$
Then
\begin{eqnarray}
\nonumber d^1_f(\Omega^n_f(\frac{g}{hf^j}))&= &
\displaystyle\sum_{i=0}^{n}d^1_f(\gfrac{g/h}{f^j,(XW)^{i+1},(YZ)^{n+1-i}})
\nonumber  \\
\label{eq:49821} &= &
\displaystyle\sum_{i=0}^{n}\gfrac{g}{hW^{i+1}Z^{n+1-i},f^j,X^{i+1},Y^{n+1-i}}
\\
\label{eq:76389} &=&
\displaystyle\sum_{i=0}^{n}
\gfrac{g(\alpha_{11}\alpha_{22}-\alpha_{12}\alpha_{21})}
{Z^{s-i},W^{t-n+i},X^{i+1},Y^{n+1-i}} \\
\nonumber &=&
g(\alpha_{11}\alpha_{22}-\alpha_{12}\alpha_{21})
\Omega^n(Z^{n+1-s}W^{n+1-t})\in E_n(Z,W).
\end{eqnarray}
The equality (\ref{eq:49821}) holds because $W$ and $Z$ are
invertible and (\ref{eq:76389}) is due to the transformation law in
Section~\ref{sec:genfrac}.
\item[Case 2. $f=Z$.] For $s > 0$,
   $\displaystyle{\Omega_Z^n(\frac{g}{h}Z^s)=0}$. For $s\leq n$, the
   elements $h$ and $Z^{n-s+1}$ form a
   system of parameters, so one may choose $t>0$ and
$\alpha_{ij}\in\kappa[Z,W]_{(Z,W)}$ such that
$$
W^t=\alpha_{21}h+\alpha_{22}Z^{n-s+1}.
$$
Then, a similar computation as in the previous case shows
\begin{eqnarray*}
d^1_Z(\Omega^n_Z(\frac{g}{h}Z^s))&=&
\displaystyle\sum_{i=0}^{n}d^1_f(\gfrac{g/h}{Z^{n+1-i-s},(XW)^{i+1},Y^{n+1-i}})\\
&=&
\displaystyle\sum_{i=0}^{n}\gfrac{g}{hW^{i+1},Z^{n+1-i-s},X^{i+1},Y^{n+1-i}}\\
&=&
\displaystyle\sum_{i=0}^{n}
\gfrac{g\alpha_{21}}{W^{t+i+1},Z^{n+1-i-s},X^{i+1},Y^{n+1-i}}\\
&=&
-g\alpha_{21}\Omega^n(Z^sW^{-t})\in E_n(Z,W).
\end{eqnarray*}
\item[Case 3. $f=W$.] Similar to Case 2.
\end{trivlist}
\end{proof}
So we have a restriction (by abusing the notation)
$$
d^1_f\colon E(f)\to E(Z,W).
$$
\begin{defn}
We define
$$
d^1\colon\underset{\scriptscriptstyle f\neq 0}{\oplus}E(f)\to E(Z,W)
$$
to be $d^1=\underset{\scriptscriptstyle f\neq 0}{\oplus}d^1_f$,
where $f$ ranges over irreducible polynomials in $(Z,W)\kappa[Z,W]$.
\end{defn}

\begin{prop}\label{prop:68500}
$d^1\circ d^0=0$.
\end{prop}
\begin{proof}
We apply the argument in the proof of \cite[Prop. 1]{hu:trpp}.
Recall that an arbitrary element in $E(0)$ can be written as a sum
of elements in the form of $\gfrac{g/h}{X^j,Y^k}$. It is enough to
show the image of such an element under $d^1 \circ d^0$ is zero in
$E(Z,W)$. We write $h=f_1\cdots f_n$ where $f_1, \dots, f_n$ are
powers of distinct irreducible polynomials. It suffices to show
\begin{equation}
\label{eq:993477}
\sum_{i=1}^n
\gfrac{g}{f_1\cdots\hat{f_i}\cdots f_n,f_i,X^j,Y^k}=0.
\end{equation}
We induct on $n$ to prove that (\ref{eq:993477}) holds for a
more general case where $f_1, \dots, f_n$ are assumed to be
products of powers of irreducible polynomials but each irreducible
factor appears in only one $f_i$.
The case $n=2$ is trivial. Assume that $n=3$. If some $f_i\not\in(Z,W)$,
(\ref{eq:993477}) clearly holds. So we assume all $f_i\in(Z,W)$.
That $f_1$ and $f_2$ are a system of parameters for
$\kappa[Z,W]_{(Z,W)}$ implies $f_3^{\ell}$ is in $(f_1,f_2)$ for some
$\ell \gg 0$. By
multiplying $g$ and $f_3$ by $f_3^{\ell-1}$ and replacing them by
the latter elements, we may assume that
$f_3=g_1 f_1+g_2 f_2$ for some $g_1,g_2\in \kappa[Z,W]_{(Z,W)}$. Then
\begin{eqnarray*}
& & \gfrac{g}{f_1 f_2,f_3,X^j,Y^k} \\
&=&\gfrac{g(g_1 f_1+g_2 f_2)}{f_1 f_2,(g_1 f_1+g_2 f_2)^2,X^j,Y^k} \\
&=&\gfrac{gg_1}{f_2,(g_1 f_1+g_2 f_2)^2,X^j,Y^k}+
\gfrac{gg_2}{f_1,(g_1 f_1+g_2 f_2)^2,X^j,Y^k} \\
&=&\gfrac{g}{f_2,g_1f_1^2,X^j,Y^k}+\gfrac{g}{f_1,g_2 f_2^2,X^j,Y^k}
\end{eqnarray*}
and
\begin{eqnarray*}
& &
\gfrac{g}{f_2 f_3,f_1,X^j,Y^k}+\gfrac{g}{f_1 f_3,f_2,X^j,Y^k}
+\gfrac{g}{f_1 f_2,f_3,X^j,Y^k} \\
&=&
\gfrac{g}{f_2^2 g_2,f_1,X^j,Y^k}+\gfrac{g}{f_1^2 g_1,f_2,X^j,Y^k}+
\gfrac{g}{f_1f_2,f_3,X^j,Y^k} \\ &= & 0.
\end{eqnarray*}
Now assume that $n>3$ and (\ref{eq:993477}) holds for numbers of $f_i$'s less

than $n$ of the general case stated above. Then
\begin{eqnarray}\label{eq:65781}
&&
\gfrac{g}{f_4\cdots f_n,(f_1 f_2 f_3),X^j,Y^k}+
\gfrac{g}{(f_1 f_2 f_3)f_5\cdots f_n,f_4,X^j,Y^k} +\cdots+\\
&&
\gfrac{g}{(f_1 f_2 f_3)f_4\cdots f_{n-1},f_n,X^j,Y^k}=0,\nonumber
\end{eqnarray}
\begin{eqnarray}\label{eq:65782}
&&\gfrac{g}{f_3(f_4\cdots f_n),(f_1 f_2),X^j,Y^k}+
\gfrac{g}{(f_1 f_2)(f_4\cdots f_n),f_3,X^j,Y^k}+\\
&&\gfrac{g}{(f_1 f_2)f_3,(f_4\cdots f_n),X^j,Y^k}=0,\nonumber
\end{eqnarray}
and
\begin{eqnarray}\label{eq:65783}
&&
\gfrac{g}{f_2(f_3\cdots f_n),f_1,X^j,Y^k}+\gfrac{g}{f_1(f_3\cdots 
f_n),f_2,X^j,Y^k}+ \\
&&
\gfrac{g}{f_1 f_2,(f_3\cdots f_n),X^j,Y^k}=0.\nonumber
\end{eqnarray}
Add identities (\ref{eq:65781}), (\ref{eq:65782}) and 
(\ref{eq:65783}), we get identity
(\ref{eq:993477}).
\end{proof}
\begin{defn}
Let $E^\bullet$ be the complex
$$
E(0)\xrightarrow{d^0}\underset{\scriptscriptstyle f\neq
0}{\oplus}E(f)\xrightarrow{d^1}E(Z,W)\to 0\to \cdots
$$
and $E^\bullet_n$
$$
E_n(0)\to\underset{\scriptscriptstyle f\neq 0}{\oplus}E_n(f)\to 
E_n(Z,W)\to 0\to\cdots
$$
be its restriction.
\end{defn}

\begin{lem}\label{89843}
$d^1$ is surjective. Let $f\in\kappa[Z,W]$ be an irreducible 
polynomial in $(Z,W)$. Then
$d^1_fE_0(f)=E_0(Z,W)$.
\end{lem}
\begin{proof}
$d^1$ is surjective, since the generators $\Omega^n(Z^sW^t)$ of
$E(Z,W)$ are in the image of $d^1$ as seen in (\ref{eq:88068}).

To prove the second assertion, we assume that $(f)\neq(W)$ to avoid
the trivial case. For any $s,t\leq 0$, we choose $\ell>0$ and
$g\in\kappa[Z,W]$ as in the proof of Lemma~\ref{prop:onto} such that
\begin{equation}\label{eq:09689}
\gfrac{g}{W^{1-t},f^\ell}=\gfrac{1}{W^{1-t},Z^{1-s}}.
\end{equation}
Since $(f) \neq (W)$, there exists $n\geq 1-s$
such that $Z^n$ is a combination of $W^{1-t}$ and $f^{\ell}$ over
$\kappa[Z,W]_{(Z,W)}$, that is,
\begin{equation}\label{eq:50832}
Z^n=\varphi W^{1-t}+\eta f^\ell
\end{equation}
for some $\varphi,\eta\in\kappa[Z,W]_{(Z,W)}$.
Using (\ref{eq:09689}) and (\ref{eq:50832}), we observe the following
$$ \begin{array}{rcl}
\gfrac{\eta g}{W^{1-t}, Z^n} &=& \gfrac{\eta g}{W^{1-t}, Z^n-\varphi
   W^{1-t}} \\ \\ & = & \gfrac{g}{W^{1-t}, f^{\ell}}
    =   \gfrac{1}{W^{1-t}, Z^{1-s}} =
    \gfrac{Z^{n+s-1}}{W^{1-t},Z^n}. \end{array}
$$
This implies
\begin{equation}\label{eq:12894}
\eta g-Z^{n+s-1}\in(W^{1-t},Z^n)
\end{equation}
in $\kappa[Z,W]_{(Z,W)}$.
The relations in (\ref{eq:50832}) and (\ref{eq:12894}) can be extended
to $S_{(X,Y,Z,W)}$. Therefore,
the second assertion follows from the computation:
$$
d^1_f(\Omega^0_f(\frac{-gZW^t}{f^\ell}))=\gfrac{-g}{W^{1-t},f^\ell,X,Y}
=\gfrac{\eta g}{Z^n,W^{1-t},X,Y}=\Omega^0(Z^sW^t).
$$
\end{proof}
\begin{lem}\label{exact1}
$E^\bullet_0$ is exact.
\end{lem}
\begin{proof}
   We only need to prove that an element of $\oplus_{f\neq 0}E_0(f)$ is
   in the image of $d^0$ if it is in the kernel of $d^1$.  Working on
   the polynomial ring $\kappa(Z)[W]$ and using Gauss lemma, one sees
   that elements in $\kappa(Z,W)$ can be written as a partial fraction
$$
\frac{g_0}{h_0}+\frac{g_1}{h_1f_1}+\cdots+\frac{g_s}{h_sf_s},
$$
where $g_i\in\kappa[Z,W]$, $h_i\in\kappa[Z]$ and $f_i\in\kappa[Z,W]$ is a
power of irreducible polynomial. This implies that, if $Z 
\not\in(f)$, elements of
$E_0(f)$ can be written as
\begin{equation}\label{eq:37862}
\Omega^0_f(\varphi Z^sf^t),
\end{equation}
where $\varphi\in\kappa[Z,W]_{(Z,W)}$, $s\in\mathbb Z$ and $t\leq-1$. Since
$$
d^0\Omega^0_0(\varphi Z^sf^t)=\Omega^0_f(\varphi Z^sf^t)+
\Omega^0_Z(\varphi Z^sf^t),
$$
and $\oplus_{f\neq 0} E_0(f)$ is generated by the image of
$\Omega^0_f$ (Definition \ref{defn:43890}), to prove the kernel of
$d^1$ contained in the image of $d^0$, we may reduce it to the case
that an element $\Psi\in E_0(Z)$ with $d^1\Psi=0$ is in the image of
$d^0$. Working on $\kappa(W)[Z]$ instead of $\kappa(Z)[W]$, we may
replace $Z$ by $W$ and choose $f$ to be $Z$ in
(\ref{eq:37862}).
Multiplying $\Psi$ by an element in $\kappa[Z,W]\setminus(Z,W)$,
we may assume
$$
\Psi=\sum_{\stackrel{\scriptstyle{s\leq 0}}{t \in \mathbb Z }}
a_{st}\Omega^0_Z(Z^sW^t)
$$
for some $a_{st} \in\kappa$ and write the map explicitly:
\begin{eqnarray*}
d^0\left(\sum_{\stackrel{\scriptstyle{s\leq 0}}{t >0}}
a_{st}\Omega^0_0(Z^sW^t) \right)
& = &
\sum_{\stackrel{\scriptstyle{s\leq 0}}{t>0}}
a_{st}\Omega^0_Z(Z^sW^t)+ \sum_{s\leq 0<t}a_{st}\Omega^0_W(Z^sW^t) \\
&=&
\sum_{\stackrel{\scriptstyle{s\leq 0}}{t>0}}a_{st}\Omega^0_Z(Z^sW^t)
= \Psi
\end{eqnarray*}
since $\Omega^0_W(Z^sW^t)=0$ for all $t >0$.
\end{proof}

Note that $E^\bullet_n$  is not exact for $n\geq 1$. For instance,
$$
d^1\Omega^n_Z(Z^nW)=0,
$$
but $\Omega^n_Z(Z^nW)$ is not in the image of $d^0$.

The maps $XW$ and $\oplus_{f\neq 0}
f^{\triangle}$ are surjective by (\ref{eq:388875}). Moreover, 
$E_0(0)$ is in the kernel of
the composition $(\oplus_{f\neq 0} f^{\triangle}) \circ d^0$ and
$\underset{\scriptscriptstyle f\neq 0}{\oplus}E_0(f)$ is in the kernel
of the composition $(X\oplus Y)\circ d^1$. We make the following
definition.
\begin{defn}
We define
$$
\pi^0\colon E(0)\to\underset{\scriptscriptstyle f\neq 0}{\oplus}E(f)
$$
and
$$
\pi^{11}\oplus\pi^{12}\colon\underset{\scriptscriptstyle f\neq 
0}{\oplus}E(f)\to E(Z,W)^2
$$
to be the maps making following diagram commutative
\begin{equation}\label{diagram}
\renewcommand\labelstyle{\scriptstyle}
\begin{diagram}
E(0)&\rTo^{\pi^0}& \underset{\scriptscriptstyle f\neq 
0}{\oplus}E(f)&\rTo^{\pi^{11}\oplus\pi^{12}}&
E(Z,W)^2 \\
\uTo^{XW} &     & \uTo_{\oplus f^{\triangle}}    &     & \uTo_{X\oplus Y}  \\
E(0)  & \rTo^{d^0}  &\underset{\scriptscriptstyle f\neq 
0}{\oplus}E(f)& \rTo^{d^1} & E(Z,W).
\end{diagram}
\end{equation}
\end{defn}

$\pi^0$ and $\pi^{11}\oplus\pi^{12}$ can be described using the maps
\begin{eqnarray*}
\pi^0_f:E(0)&\to& E(f),\\
\pi^{11}_f:E(f)&\to&E(Z,W),\\
\pi^{12}_f:E(f)&\to&E(Z,W),
\end{eqnarray*}
where
$$
\begin{cases}
\pi^0_Z=d^0_Z\circ\tfrac{1}{Z},\\
\pi^0_W=d^0_W\circ\tfrac{1}{W},\\
\pi^0_f=d^0_f,
\end{cases}
\begin{cases}
\pi^{11}_Z=d^1_Z\circ\tfrac{Z}{W},\\
\pi^{11}_W=d^1_W,\\
\pi^{11}_f=d^1_f\circ\tfrac{1}{W},
\end{cases}
\begin{cases}
\pi^{12}_Z=d^1_Z,\\
\pi^{12}_W=d^1_W\circ\tfrac{W}{Z},\\
\pi^{12}_f=d^1_f\circ\tfrac{1}{Z},
\end{cases}$$
for $Z,W\not\in(f)$. We have
\begin{eqnarray*}
\pi^{11}=\underset{\scriptscriptstyle f\neq 0}{\oplus}\pi^{11}_f
&\text{ and } &
\pi^{12}=\underset{\scriptscriptstyle f\neq 0}{\oplus}\pi^{12}_f.
\end{eqnarray*}
The product $\prod_{f\neq 0}\pi^0_f\colon E(0)\to \prod_{f\neq 
0}E(f)$ has image in $\oplus_{f\neq 0}E(f)$ and
equals $\pi^0$.

\begin{defn}
We define
\begin{equation}\label{eq:79989}
E(0)
\xrightarrow{\delta^0}\underset{\scriptscriptstyle \q\neq(Z,W)}{\oplus}E(\q)
\xrightarrow{\delta^1}\underset{\scriptscriptstyle \q\neq(0)}{\oplus}E(\q)
\xrightarrow{\delta^2}E(Z,W)^2
\end{equation}
to be the total complex associated to the double complex
(\ref{diagram}) with a negative sign on $\oplus f^{\triangle}$.
For $n\geq 3$, we define
$$
\delta^n\colon E(Z,W)^2\to E(Z,W)^2
$$
to be the map
$$
\delta^n(\Psi_1\oplus\Psi_2)=
\begin{cases}
(W\Psi_1-Z\Psi_2)\oplus(-Y\Psi_1+X\Psi_2),&\text{if $n$ is odd;}\\
(X\Psi_1+Z\Psi_2)\oplus(Y\Psi_1+W\Psi_2),&\text{if $n$ is even.}
\end{cases}
$$
\end{defn}

\begin{thm}\label{thm:main}
\begin{small}\begin{equation}
\label{eq:main}
E(0)
\xrightarrow{\delta^0}\underset{\scriptscriptstyle\q\neq\m}{\oplus}E(\q)
\xrightarrow{\delta^1}\underset{\scriptscriptstyle\q\neq(0)}{\oplus}E(\q)
\xrightarrow{\delta^2}E(Z,W)^2\xrightarrow{\delta^3}E(Z,W)^2\xrightarrow{\delta^4}E(Z,W)^2\cdots
\end{equation}\end{small}
is a minimal injective resolution of $A/\p$.
\end{thm}
\begin{proof}
$\pi^0\circ(XW)- (\oplus f^{\triangle})
  \circ d^0=0$, by the definition of $\pi^0$. This identity, together with
Proposition~\ref{prop:68500} (that is, $d^1\circ d^0=0$), implies that
$\delta^1\circ\delta^0=0$. The multiplication by $XW$ on $E(0)$ is 
surjective. Therefore to
show $\delta^2\circ\delta^1=0$, we only need to check whether the 
image of an element of
$\oplus_{f\neq 0}E(f)$ vanishes. This is easy, since
$(X\oplus Y)\circ d^1-(\pi^{11}\oplus\pi^{12})\circ\mu=0$ by the definition of
$\pi^{11}\oplus\pi^{12}$. The map $\mu$ is surjective. Therefore to show
$\delta^3\circ\delta^2=0$, we only need to compute the image of an element of
$E(Z,W)$. This is also easy, since elements of $E(Z,W)$ are 
annihilated by $XW-YZ$.
It is straightforward to show that $\delta^{n+1}\circ\delta^n=0$ for 
$n\geq 3$. We conclude
that (\ref{eq:main}) is a complex.

The maps $XW$ and $-(\oplus f^{\triangle})$ are surjective. Chasing 
diagram (\ref{diagram}),
it is easy to see that (\ref{eq:79989}) is exact because
$$
E_0(0)\to\underset{\scriptscriptstyle f\neq 0}{\oplus}E_0(f)\to E_0(Z,W)\to 0
$$
is exact. The complex
$$
\cdots A^2
\xrightarrow{\tiny\begin{pmatrix}X&Y\\Z&W\end{pmatrix}}A^2
\xrightarrow{\tiny\begin{pmatrix}W&-Y\\-Z&X\end{pmatrix}}A^2
\xrightarrow{\tiny\begin{pmatrix}X&Y\\Z&W\end{pmatrix}}A^2
\xrightarrow{\tiny\begin{pmatrix}W&-Y\\-Z&X\end{pmatrix}}A^2
\xrightarrow{\tiny\begin{pmatrix}X&Y\end{pmatrix}}A
$$
is exact. Apply the functor $\Hom_A(-,E(Z,W))$, we get the exact sequence
$$
E(Z,W)\xrightarrow{X\oplus Y}E(Z,W)^2
\xrightarrow{\delta^3}E(Z,W)^2\xrightarrow{\delta^4}E(Z,W)^2\xrightarrow{\delta^5}E(Z,W)^2\to\cdots.
$$
Hence
\begin{equation}\label{eq:90688}
\underset{\scriptscriptstyle \q\neq(0)}{\oplus}E(\q)
\xrightarrow{\delta^2}E(Z,W)^2\xrightarrow{\delta^3}E(Z,W)^2\xrightarrow{\delta^4}E(Z,W)^2\xrightarrow{\delta^5}E(Z,W)^2\to\cdots
\end{equation}
is also exact. Combining (\ref{eq:79989}) and (\ref{eq:90688}), we 
conclude that
(\ref{eq:main}) is exact.

The kernel of $\delta^0$ equals the kernel of $d^0$ restricting to 
$E_0(0)$. Let
$$
\Psi=\gfrac{g/h}{XW,YZ}
$$
be an element of $E_0(0)$ in the kernel of $d^0$, where $g\in\kappa[Z,W]$ and
$0\neq h\in\kappa[Z,W]$. Then $g\in hZW\kappa[Z,W]_{(f)}$ for any 
irreducible polynomial
$f\in(Z,W)$. Therefore $g=\varphi hZW$ for some 
$\varphi\in\kappa[Z,W]_{(Z,W)}$ and
\begin{equation}\label{eq:79622}
\Psi=\gfrac{\varphi}{X,Y}.
\end{equation}
All elements of the above form is in the kernel of $d^0$. These 
elements form a module isomorphic to
$A/\p$. Therefore (\ref{eq:main}) is an injective resolution of $A/\p$.

Let $\Psi=\sum_{i=0}^n \Omega^i_0(\varphi_i)$ be a non-zero element of
$E(0)$ with $\varphi_n \neq0$ and $\varphi_n= g/h$ for some
$g, h \in \kappa[Z,W]$. By the structure of $E(f)$
discussed in Example \ref{ex:E(f)},
$$
hWZ(XW)^n \Psi = hWZ \Omega_0^0(\varphi_n) = \gfrac{g}{X,Y}
$$
which is a non-zero element in $E_0(0)$ of the form (\ref{eq:79622}),
so it is in the image of $A/ \p$. Hence $E(0)$ is an injective hull of $A/\p$.

Every non-zero element of $\oplus_{\q\neq(Z,W)}E(\q)$ multiplied by a
suitable element of $\kappa[Z,W]$ becomes a non-zero element in the
summand $E(0)$ of the form $\sum \Omega_0^i(g_i)$, $g_i \in
\kappa[Z,W]$. This element is in the image of $\delta^0$.

Every non-zero element of $\oplus_{f\neq 0}E(f)$ multiplied by
suitable powers of irreducible polynomials becomes a non-zero element
$\Psi_1=\sum \Omega^i_f(h_i/f^{r_i})\in E(f)$ for some non-zero $f,
h_i \in \kappa[Z,W]$, and $n$, $r_i \in \mathbb N$.  By definition,
$$\sum \Omega^i_f(\frac{h_{i}}{f^{r_i}})=d^0_f \sum 
\Omega_0^i(\frac{h_{i}}{f^{r_i}})=\pi^0 \sum\Omega_0^i(\frac{h_{i}}{f^{r_i}})$$
if $Z,W\not\in(f)$ and $$\sum\Omega^i_f(\frac{h_{i}}{f^{r_i}})=
d^0_f \circ \frac{1}{f}( f \sum\Omega_0^i(\frac{h_{i}}{f^{r_i}})) 
=\pi^0 \sum f\Omega_0^i(\frac{h_{i}}{f^{r_i}})$$
if $f=Z$ or $W$. So $(\Psi_1,0)$ is in the image of $\delta^1$. Every non-zero
element of $E(Z,W)$ multiplied by suitable powers of $X$ and $Y$ 
becomes a non-zero
element $\Psi_2\in E_0(Z,W)$ (see Example \ref{ex:E(Z,W)}). By Lemma
\ref{89843} and the definition of $f^{\triangle}$, $(0, \Psi_2)$ is in
the image of $\delta ^1$.  Now for a general case, multiplied by
suitable powers of irreducible polynomials in $\kappa[Z,W]$, and those
of $X$ and $Y$, a non-zero element of $(\oplus_{f\neq 0}E(f)) \oplus
E(Z,W)$ becomes a non-zero element $(\Psi_1, \Psi_2)$ with $\Psi_1,\Psi_2$ as
described in the above and therefore, it is in the image of $\delta^1$.

Every non-zero element of $E(Z,W)^2$ multiplied by suitable powers of 
$X$ and $Y$
becomes a non-zero element of $E_0(Z,W)^2$.  Multiplied again by 
suitable powers of $Z$ and
$W$, this element becomes a non-zero element of the form
\begin{eqnarray*}
\alpha\Omega^0(1)\oplus\beta\Omega^0(1) && (\alpha,\beta\in\kappa),
\end{eqnarray*}
which is in the image of $\delta^n$ ($n\geq 2$), since
\begin{eqnarray*}
&&
\begin{cases}
\pi^{11}(\alpha\Omega^0_W(1)\oplus-\beta\Omega^0_Z(1))=\alpha\Omega^0(1)\\
\pi^{12}(\alpha\Omega^0_W(1)\oplus-\beta\Omega^0_Z(1))=\beta\Omega^0(1)
\end{cases}\\
&&\begin{cases}
\alpha X\Omega^1(W)+\beta Z\Omega^0(W^{-1})=\alpha\Omega^0(1)\\
\alpha Y\Omega^1(W)+\beta W\Omega^0(W^{-1})=\beta\Omega^0(1)
\end{cases}
\\
&&\begin{cases}
\alpha W\Omega^0(W^{-1})-\beta Z\Omega^1(W)=\alpha\Omega^0(1)\\
-\alpha Y\Omega^0(W^{-1})+\beta X\Omega^1(W)=\beta\Omega^0(1).
\end{cases}
\end{eqnarray*}
Therefore the resolution is minimal.
\end{proof}
\begin{cor}
The Bass numbers of $A/\p$ are as follows. Let $f$ be an irreducible 
polynomial contained in $(Z,W)$.
\begin{eqnarray*}
\mu_i((X,Y,Z,W),A/\p)&=&
\begin{cases}0,&\text{if $i<2$;}\\1,&\text{if $i=2$;}\\2, &\text{if 
$i>2$.}\end{cases}\\
\mu_i((X,Y),A/\p)&=&
\begin{cases}1,&\text{if $i<2$;}\\0, &\text{if $i\geq 2$.}\end{cases}\\
\mu_i((X,Y,f),A/\p)&=&
\begin{cases}0,&\text{if $i=0$;}\\1,&\text{if $i=1$ or $2$;}\\0, 
&\text{if $i>2$.}\end{cases}
\end{eqnarray*}
All other Bass numbers of $A/\p$ are zero.
\end{cor}

Minimal injective resolutions of $A/\p$ are eventually periodic. 
K.~Yanagawa informs
us that this is true in a general setting:  Over a local ring which 
is a hypersurface
with an isolated singularity, minimal injective resolutions of any 
finitely generated
modules are eventually periodic. His proof uses Matlis duality and a result of
Eisenbud \cite{eis:haci}


\section{Applications}\label{sec:cohomology}


\subsection{Local Cohomology}
We compute the local cohomology module $\mH^i_I(A/\p)$ of $A/\p$
supported at an ideal $I$ of $A$. Recall that the injective resolution
(\ref{eq:main}) of $A/\p$ is built up by injective hulls $E(\q)$ of
modules $A/(\q,X,Y)$, where $\q$ is a prime ideal of $\kappa[Z,W]$.
Since elements in $E(\q)$ are annihilated by powers of $X$ and $Y$,
the functors $\Gamma_I$ and $\Gamma_{I+(X,Y)}$ have the same effect on
the complex (\ref{eq:main}). Hence
$\mH^i_I(A/\p)=\mH^i_{I+(X,Y)}(A/\p)$ for all $i$ and, to compute the
local cohomology modules, we may assume that $I=(I_0,X,Y)$ for some
ideal $I_0$ of $\kappa[Z,W]$. If $I_0\subset\q$, then $E(\q)$ being
$\q$-torsion is also $I$-torsion. If $I_0\not\subset\q$, there is an
element $a\in I_0\setminus\q$. The only element of $E(\q)$ annihilated
by powers of $a$ is zero, so $E(\q)$ is $I$-torsion
free in this case. Therefore applying the $I$-torsion functor
$\Gamma_I(-)$ simply means taking away those $E(\q)$ with
$I_0\not\subset\q$ from the complex (\ref{eq:main}).

If $\height(\kappa[Z,W]\cap I)=2$, $(X,Y,Z,W)$ is the only prime
containing $I$. Apply the functor $\Gamma_I(-)$ to the injective
resolution (\ref{eq:main}) of $A/\p$, we get the complex
$$
0\to 0\to E(Z,W)\xrightarrow{X\oplus Y}E(Z,W)^2
\xrightarrow{\delta^3}E(Z,W)^2\xrightarrow{\delta^4}E(Z,W)^2\to\cdots
$$
whose only non-trivial cohomology is $E_0(Z,W)$. Therefore
\begin{equation}
\mH^i_I(\frac{A}{\p})=
\begin{cases}E_0(Z,W), & \text{if $i=2$;}\\0, & \text{if $i\neq 2$.}\end{cases}
\end{equation}
As a $\kappa$-vector space, $\mH^2_I(A/\p)$ has a basis consisting of 
$\Omega^0(Z^sW^t)$, where
$s,t\leq 0$.

A local cohomology module is a direct limit of extension modules. 
Using the injective resolution
(\ref{eq:main}), we can see clearly the behavior of the limit
$$
\lim_{n\to\infty}\Ext_A^2(\frac{A}{(X,Y,Z,W)^n},\frac{A}{\p})=\mH^2_{(X,Y,Z,W)}(\frac{A}{\p}).
$$
Apply the functor $\Hom_A(A/(X,Y,Z,W)^n,-)$ to  (\ref{eq:main}), we get
\begin{small}$$
   0\to 
0\to\Hom_A(\frac{A}{(X,Y,Z,W)^n},E(Z,W))\xrightarrow{\tiny\begin{pmatrix}X\\Y\end{pmatrix}}
   \Hom_A(\frac{A}{(X,Y,Z,W)^n},E(Z,W)^2)\to\cdots.
$$\end{small}
$\Ext_A^2(A/(X,Y,Z,W)^n,A/\p)$ is isomorphic to the submodule of 
$E(Z,W)$ consisting of those
elements annihilated by $X$, $Y$ and $(X,Y,Z,W)^n$. As a $\kappa$-vector space,
$\Ext_A^2(A/(X,Y,Z,W)^n,A/\p)$ has a basis consisting of
$\Omega^0(Z^sW^t)$, where $s,t\leq 0$ satisfy $s+t+n>0$. In particular,
$$
\dim_\kappa\Ext^2(\frac{A}{(X,Y,Z,W)^n},\frac{A}{\p})=\frac{n(n+1)}{2}.
$$
As $n$ increasing, the set
$\{\Omega^0(Z^sW^t)|\text{ $s,t\leq 0$},\text{ $s+t+n>0$}\}$ becomes 
larger and closer to the
basis $\{\Omega^0(Z^sW^t)|\text{ $s,t\leq 0$}\}$ of $\mH^2_{(X,Y,Z,W)}(A/\p)$.

If $\height(\kappa[Z,W]\cap I)=0$, then $I=(X,Y)$. The functor 
$\Gamma_I(-)$ does not
change the complex (\ref{eq:main}). Therefore
\begin{equation}
\mH^i_{(X,Y)}(\frac{A}{\p})=
\begin{cases}A/\p, & \text{if $i=0$;}\\0, & \text{if $i\neq 0$.}\end{cases}
\end{equation}

Now we look at the case $\height(\kappa[Z,W]\cap I)=1$. Applying 
$\Gamma_I(-)$ to
(\ref{eq:main}), we get a complex quasi-isomorphic to
$$
  0\to
\underset{\scriptscriptstyle I\subset(f,X,Y)}{\oplus}
E_0(f)\to E_0(Z,W)\to 0\to \cdots.
$$
By Lemma~\ref{89843}, the non-trivial map in the above complex is 
surjective. Therefore
\begin{equation}
\mH^i_I(\frac{A}{\p})=
\begin{cases}\text{kernel of $\oplus_{I\subset(f,X,Y)}E_0(f)\to E_0(Z,W)$}, &
\text{if
$i=1$;}\\0, & \text{if $i\neq 1$.}\end{cases}
\end{equation}
For instance, if
$I=(Z,X,Y)$, the above complex becomes
$$
\cdots\to 0\to E_0(Z)\to E_0(Z,W)\to 0\to \cdots.
$$
As a $\kappa[W]_{(W)}$-module, $\mH^1_{(Z,X,Y)}(A/\p)$ is generated freely by
$\Omega^0_Z(Z^sW)$, where $s\leq 0$.

\subsection{Normal Module}\label{sebsec:86}
We would like to make explicit the canonical isomorphism
\begin{equation}\label{eq:894421}
\Hom_A(\p/\p^2,A/\p)\to\Ext^1_A(A/\p,A/\p)
\end{equation}
in terms of the injective resolution (\ref{eq:main}) of $A/\p$. Note 
that the canonical map
$$
\Hom_A(\p/\p^2,A/\p)\to\Hom_A(\p,A/\p)
$$
is an isomorphism. We describe an isomorphism between
$\Ext^1_A(A/\p,A/\p)$ and $\Hom_A(\p,A/\p)$ to establish (\ref{eq:894421}).

We compute $\Ext^1_A(A/\p,A/\p)$ by applying the functor 
$\Hom_A(A/\p,-)$ to (\ref{eq:main}).
By Proposition ~\ref{prop:086661}, $\Ext^1_A(A/\p,A/\p)$ is the 
cohomology of the complex
$$
E_0(0)\xrightarrow{\delta^0}\underset{\scriptscriptstyle\q\neq(Z,W)}{\oplus}E_0(\q)
\xrightarrow{\delta^1}\underset{\scriptscriptstyle\q\neq(0)}{\oplus}E_0(\q).
$$
By Lemma~\ref{exact1}, it is also the kernel of
$$
E_0(0)\xrightarrow{\pi^0}\underset{\scriptscriptstyle f\neq 0}{\oplus}E_0(f).
$$
Explicitly,
$$
\Ext^1_A(A/\p,A/\p)\simeq\{\Omega^0_0(gZ^2W^2)|g\in\kappa[Z,W]_{(Z,W)}\}\subset 
E_0(0).
$$
Consider the diagram
$$
\renewcommand\labelstyle{\scriptstyle}
\begin{diagram}
            &           &           &          & \Hom(\p,A/\p)\\
            &           &           &          & \dTo \\
E_0(0)& \rTo  & E(0)   & \rTo & \Hom(\p,E(0))\\
\dTo   &           & \dTo  &         & \dTo \\
\underset{\scriptscriptstyle\q\neq(Z,W)}{\oplus}E_0(\q)& \rTo
&\underset{\scriptscriptstyle\q\neq(Z,W)}{\oplus}E(\q)&\rTo&
\Hom(\p,\underset{\scriptscriptstyle\q\neq(Z,W)}{\oplus}E(\q))\\
\dTo   &           & \vdots  &         & \vdots\\
\underset{\scriptscriptstyle\q\neq(0)}{\oplus}E_0(\q)& && &\\
\vdots   &           &     &         &
\end{diagram}
$$
obtained by applying the $\Hom$ functors on the short exact sequence
$$
0\to\p\to A\to A/\p\to 0
$$
and using (\ref{eq:main}) to establish the vertical maps.
Chasing the above diagram, we get an isomorphism
$$
\Ext^1_A(A/\p,A/\p)\to\Hom_A(\p,A/\p),
$$
which maps
$\Omega^0_0(gZ^2W^2)$ to the $A$-linear map $\p\to A/\p$ determined by
\begin{eqnarray*}
X\mapsto gZ &\text{and} &Y\mapsto gW
\end{eqnarray*}
for $g\in\kappa[Z,W]_{(Z,W)}$.


\subsection{Yoneda algebra}\label{subsec:Yoneda}


First we compute $\Ext^*_A(A/\p,A/\p)=\sum_{i=0}^\infty\Ext^i_A(A/\p,A/\p)$. In
Subsection~\ref{sebsec:86}, we have seen that
$$
\Ext^1_A(A/\p,A/\p)\simeq\{\Omega^0_0(gZ^2W^2)|g\in\kappa[Z,W]_{(Z,W)}\}.
$$
Let
\begin{eqnarray*}
e_0 &=& \Omega^0_0(ZW)\in E(0),\\
e_1 &=& \Omega^0_0(Z^2W^2)\in\underset{\scriptscriptstyle\q\neq\m}{\oplus}E(\q),\\
e_2 &=& \Omega^0_W(Z)\in\underset{\scriptscriptstyle\q\neq 0}{\oplus}E(\q),\\
e_{2n} &=& 0\oplus\Omega^0(1) \text{ in the $2n$-th term $E(Z,W)^2$ 
of (\ref{eq:main}),}
\end{eqnarray*}
where $n>1$, be (the equivalence classes of) the cycles of the complex
(\ref{eq:main}).
It should be pointed out that all the above $e_i$ represent
non-trivial cohomology classes in $\Ext^i_A(A/\p, A/\p)$.
For each $e_j$, we define
$$
\iota_j\colon A/\p\to\text{ the $j$-th term of (\ref{eq:main})}
$$
to be the map sending $1$ to $e_j$, in which $\iota_0$ is the embedding making
(\ref{eq:main}) in Theorem~\ref{thm:main} an injective resolution of $A/\p$.

The following Lemma~\ref{lem:588866} describes $\Ext^*(A/\p,A/\p)$ as
an $A$-module using independent generators $e_i$. Later in
Proposition~\ref{prop:Yoneda}, we will present $\Ext^*(A/\p,A/\p)$ as an $A$-algebra.
\begin{lem}\label{lem:588866}
As an $A$-module, the Yoneda algebra $\Ext^*_A(A/\p,A/\p)$ is generated by 
$e_0,e_1,e_2,e_4,e_6,\cdots$. The annihilators of $e_0$ and $e_1$ are $\p$; for $i>0$, the
annihilator of $e_{2i}$ is $\p+AZ+AW$.
\end{lem}
\begin{proof}
It is clear that $\Ext^0_A(A/\p,A/\p)$ is generated by $e_0$, whose annihilator is $\p$.
The module $\Ext^1_A(A/\p,A/\p)$ has been treated in Subsection~\ref{sebsec:86}.

As seen in Subsection~\ref{sebsec:86}, for all $i \geq 2$,
$\Ext^i_A(A/\p,A/\p)$ is a cohomology module of the complex
$$
E_0(0)\xrightarrow{\pi^0}\underset{\scriptscriptstyle f\neq 0}{\oplus}E_0(f)
\xrightarrow{\pi^{11}\oplus\pi^{12}}E_0(Z,W)^2\xrightarrow{\delta^3}E_0(Z,W)^2\to\cdots.
$$
For $n\geq 2$, $\Ext^{2n+1}_A(A/\p,A/\p)=0$, since the complex
$$
E_0(Z,W)^2\xrightarrow{\tiny\begin{pmatrix}0&Z\\0&W\end{pmatrix}}E_0(Z,W)^2
\xrightarrow{\tiny\begin{pmatrix}W&-Z\\0&0\end{pmatrix}}E_0(Z,W)^2
$$
is exact. The exactness of the above sequence means that the kernel of
$\tiny\begin{pmatrix}W&-Z\\0&0\end{pmatrix}$ consists of elements of 
the form $Z\Psi\oplus
W\Psi$, where $\Psi\in E_0(Z,W)$. By Lemma~\ref{89843}, these 
elements are in the image of
$W(\pi^{11}_Z\oplus\pi^{12}_Z)=(Zd^1_Z)\oplus(Wd^1_Z)$. Therefore
$$
E_0(Z)\xrightarrow{\pi^{11}_Z\oplus\pi^{12}_Z}E_0(Z,W)^2
\xrightarrow{\tiny\begin{pmatrix}W&-Z\\0&0\end{pmatrix}}E_0(Z,W)^2
$$
is exact and $\Ext^3_A(A/\p,A/\p)=0$ as well.

The cohomology of the complex
$$
E_0(Z,W)^2
\xrightarrow{\tiny\begin{pmatrix}W&-Z\\0&0\end{pmatrix}}E_0(Z,W)^2
\xrightarrow{\tiny\begin{pmatrix}0&Z\\0&W\end{pmatrix}}E_0(Z,W)^2
$$
is generated by the element in $\Ext^{2n}_A(A/\p,A/\p)$ represented by
$0\oplus\Omega^0(1)$. Therefore $\Ext^{2n}_A(A/\p,A/\p)=Ae_{2n}$
for $n\geq 2$. The element $e_{2n}$ is non-zero and annihilated by $Z$, $W$ and $\p$.
Its annihilator is hence $\p+AZ+AW$.

As seen in the proof of Lemma~\ref{exact1}, if $Z,W\not\in(f)$, elements of
$E_0(f)$ can be written as
$$
\Omega^0_f(\varphi Z^sf^t),
$$
where $\varphi\in\kappa[Z,W]_{(Z,W)}$, $s\in\mathbb Z$ and $t\leq-1$. Since
$$
\Omega^0_f(\varphi Z^sf^t)=\pi^0\Omega^0_0(\varphi Z^sf^t)-
\Omega^0_Z(\varphi Z^{s-1}f^t)-\Omega^0_W(\varphi Z^sf^tW^{-1}),
$$
to compute $\Ext^2_A(A/\p,A/\p)$, we may restrict $\pi^{11}\oplus\pi^{12}$
to $E_0(Z)\oplus E_0(W)$. Multiplied by an element in 
$\kappa[Z,W]\setminus(Z,W)$,
an element in $E_0(Z)\oplus E_0(W)$ can be written as the form
$$
\sum_{\stackrel{\scriptstyle{s\leq 0}}{t\in\mathbb
     Z}}a_{st}\Omega^0_Z(Z^sW^t) +
\sum_{\stackrel{\scriptstyle{s\in\mathbb Z}}{t\leq 0}}b_{st}\Omega^0_W(Z^sW^t).
$$
Since
$$
\Omega^0_Z(Z^sW^t)=\pi^0\Omega^0_0(Z^{s+1}W^t)-\Omega^0_W(Z^{s+1}W^{t-1})
$$
and
\begin{eqnarray*}
\Omega^0_W(Z^sW^t)=\pi^0\Omega^0_0(Z^sW^{t+1}) & \text{ for } & s>1,
\end{eqnarray*}
to compute $\Ext^2_A(A/\p,A/\p)$, we may work on elements of the form
$$
\sum_{\stackrel{\scriptstyle{s\leq 1}}{t\leq 0}}b_{st}\Omega^0_W(Z^sW^t)
$$
and assume that it is in the kernel of
$\pi^{11}\oplus\pi^{12}$.
$$\pi^{11} \left( \sum_{\stackrel{\scriptstyle{s\leq 1}}{t\leq
     0}}b_{st}\Omega^0_W(Z^sW^t) \right)=0
$$
implies $b_{st}=0$ for all $s,t \leq 0$. Furthermore,
$$\pi^{12}\left(\sum_{t\leq 0}b_{1t}\Omega_W^0(ZW^t)\right)=0$$
implies $b_{1t}=0$ for all $t\leq -1$. Therefore, $b_{10}$ is the only
possible non-zero coefficient and $\Ext^2_A(A/\p,A/\p)=Ae_2$.
Clearly, $e_2$ is annihilated by $W$ and $\p$. 
Since $Ze_2=\pi^0\Omega^0_0(Z^2W)$, it is also annihilated by $Z$.
Finally, $e_2$ is non-zero, so its annihilator is $\p+AW+AZ$.
\end{proof}

Now we compute the Yoneda pairing
$$
\Ext^i_A(A/\p,A/\p)\times\Ext^j_A(A/\p,A/\p)\to\Ext^{i+j}_A(A/\p,A/\p).
$$
Since the pairing is $A$-bilinear, we only need to compute $e_i\times e_j$.

\begin{lem}
$$
e_i\times e_j=
\begin{cases}
e_i, &\text{if $j=0$;}\\
e_j, &\text{if $i=0$;}\\
0, &\text{if $i=1$ or $j=1$, but $ij\neq 0$;}\\
-e_{i+j}, &\text{if $ij\neq 0$ and $i,j$ are both even.}
\end{cases}
$$
\end{lem}
\begin{proof}
To compute $e_2\times e_2$, we need to construct a commutative diagram
\begin{equation}\label{eq:07980}
\renewcommand\labelstyle{\scriptstyle}
\begin{diagram}
A/\p & \rTo^{\iota_0} & E(0) &\rTo^{\delta^0} 
&\underset{\scriptscriptstyle\q\neq\m}{\oplus}E(\q)
&\rTo^{\delta^1} &\underset{\scriptscriptstyle\q\neq(0)}{\oplus}E(\q)\\
&\rdTo_{\iota_2} &\dTo & &\dTo &&\dTo\\
&&\underset{\scriptscriptstyle\q\neq(0)}{\oplus}E(\q)
&\rTo_{\delta^2} & E(Z,W)^2&\rTo_{\delta^3} &E(Z,W)^2.
\end{diagram}
\end{equation}
It is straightforward to check that the diagrams
$$
\renewcommand\labelstyle{\scriptstyle}
\begin{diagram}
A/\p & \rTo^{\iota_0} & E(0) &\rTo^{d^0} 
&(\underset{\scriptscriptstyle f\neq 0,Z,W}{\oplus}E(f))\oplus E(Z)\oplus E(W)\\
&\rdTo_{\iota_2} &\dTo_{d^0_W\frac{1}{W}} & &\dTo_{\M_{23}} \\      
&&E(W)&\rTo_{\tiny\begin{pmatrix}d^1_W\\d^1_W\frac{W}{Z}\end{pmatrix}} & E(Z,W)^2
\end{diagram}
$$
and
$$
\renewcommand\labelstyle{\scriptstyle}
\begin{diagram}
(\underset{\scriptscriptstyle f\neq 0,Z,W}{\oplus}E(f))\oplus E(Z)\oplus E(W) &
\rTo^{\M_{43}} &
E(Z,W)\oplus(\underset{\scriptscriptstyle f\neq 0,Z,W}{\oplus}E(f))\oplus E(Z)\oplus E(W)\\
\dTo^{\M_{23}} &&\dTo_{\M_{24}}\\      
E(Z,W)^2
&\rTo_{\tiny\begin{pmatrix}W&-Z\\-Y&X\end{pmatrix}} &
E(Z,W)^2
\end{diagram}
$$
are commutative, where
\begin{eqnarray*}
\M_{23}&=&\begin{pmatrix}-\underset{\scriptscriptstyle f\neq 0,Z,W}{\oplus}d^1_f\frac{1}{W}&-d^1_Z\frac{1}{W}&0\\
0&0&d^1_W\frac{1}{Z}\end{pmatrix},\\
\M_{24}&=&\begin{pmatrix}-1&0&0&0\\0&-\underset{\scriptscriptstyle f\neq 0,Z,W}{\oplus}d^1_f\frac{1}{ZW}&-d^1_Z\frac{1}{W}&-d^1_W\frac{1}{Z}\end{pmatrix},\\
\M_{43}&=&\begin{pmatrix}\underset{\scriptscriptstyle f\neq 0,Z,W}{\oplus}d^1_f&d^1_Z&d^1_W\\
-\underset{\scriptscriptstyle f\neq 0,Z,W}{\oplus}f^\Delta&0&0\\0&-Z^\Delta&0\\0&0&-W^\Delta\end{pmatrix}.
\end{eqnarray*}
We define the vertical maps in (\ref{eq:07980}) to be those in the above and zero maps if a
certain component is not included above. The product $e_2\times e_2$ is the image of $e_2$ under
the map
$\M_{24}$, which equals
$-e_4$. 

We use the same method to compute other $e_i\times e_j$. For $i>1$, the diagram
$$
\begin{small}\renewcommand\labelstyle{\scriptstyle}
\begin{diagram}
A/\p & \rTo^{\iota_0} & E(0) & \rTo^{d^0} & 
\underset{\scriptscriptstyle f\neq 0}{\oplus}E(f) &
                                            \rTo^{\M_{43}} & \oplus_{\q\neq (0)}E(\q) \\
& \rdTo_{\iota_{2i}} & 
\dTo_{\tiny\begin{pmatrix}0\\d^1_Wd^0_W\frac{1}{ZW}\end{pmatrix}}  & 
&    \dTo_{\M_{23}} & & \dTo_{\M_{24}}\\
& & E(Z,W)^2 & \rTo_{\tiny\begin{pmatrix}X&Z\\Y&W\end{pmatrix}}  & E(Z,W)^2
                 & \rTo_{\tiny\begin{pmatrix}W&-Z\\-Y&X\end{pmatrix}} & E(Z,W)^2
\end{diagram}\end{small}
$$
commutes. Therefore $e_2\times e_{2i}=-e_{2i+2}$ for $i>1$. The diagram
$$
\begin{small}\renewcommand\labelstyle{\scriptstyle}
\begin{diagram}
\oplus_{\q\neq (0)}E(\q)&\rTo^{\delta^2} &E(Z,W)^2
&\rTo^{\tiny\begin{pmatrix}W&-Z\\-Y&X\end{pmatrix}} &
E(Z,W)^2&\cdots\\
\dTo^{\M_{24}}&
&\dTo_{\tiny\begin{pmatrix}-1&0\\0&-1\end{pmatrix}}&&\dTo_{\tiny\begin{pmatrix}-1&0\\0&-1\end{pmatrix}}&\\
E(Z,W)^2&\rTo_{\tiny\begin{pmatrix}X&Z\\Y&W\end{pmatrix}} &E(Z,W)^2
&\rTo_{\tiny\begin{pmatrix}W&-Z\\-Y&X\end{pmatrix}} &
E(Z,W)^2&\cdots
\end{diagram}\end{small}
$$
commutes. Therefore $e_{2j}\times e_{2i}=-e_{2i+2j}$ for $i\geq 1$ and $j>1$.
The diagram
$$
\renewcommand\labelstyle{\scriptstyle}
\begin{diagram}
A/\p & \rTo^{\iota_0}  & E(0)   & \rTo^{d^0} &
(\underset{\scriptscriptstyle f\neq 0,Z,W}{\oplus}E(f))\oplus E(Z)\oplus E(W)\\
& \rdTo_{\iota_1}  &\dTo_{ZW}  &         &
\dTo_{\tiny\begin{pmatrix}ZW&0&0\\0&W&0\\0&0&Z\end{pmatrix}} \\    &
&E(0)&\rTo_{\pi^0}&
(\underset{\scriptscriptstyle f\neq 0,Z,W}{\oplus}E(f))\oplus E(Z)\oplus E(W)
\end{diagram}
$$
commutes and $\iota_1$ has image in $E(0)$. Therefore $e_1\times e_1=0$.

$\Ext^{2i+1}_A(A/\p,A/\p)=0$ for $i\geq 1$. Therefore $e_1\times 
e_{2i}=e_{2i}\times e_1=0$
for $i\geq 1$. It is easy to see that $e_0\times e_i=e_i\times e_0=e_i$.
\end{proof}

\begin{prop}\label{prop:Yoneda}
The Yoneda algebra $\Ext^*_A(A/\p,A/\p)$ is isomorphic to the
polynomial ring $A/\p[U,V]$ modulo the ideal generated by
$ZV$, $WV$, $U^2$ and $UV$. 
\end{prop}
\begin{proof}
All $e_i$ are annihilated by $\p$. So there is an $A$-algebra homomorphism
$$
A/\p[U,V]\to\Ext^*_A(A/\p,A/\p)
$$ 
given by $1\mapsto e_0$, $U\mapsto e_1$ and $V\mapsto e_2$.
Since $\Ext^*_A(A/\p,A/\p)$ is generated by $e_i$ and
$(-1)^{n+1}e_{2n}=e_2^n$ (the Yoneda product of $n$ copies of $e_2$), the homomorphism is 
surjective. By Lemma~\ref{lem:588866}, the kernel 
of the homomorphism is generated by $ZV$, $WV$, $U^2$ and 
$UV$.
\end{proof}

\begin{cor}
The Yoneda algebra $\Ext^*_A(A/\p,A/\p)$ is commutative and finitely generated.
\end{cor}


\subsection{Dutta, Hochster and McLaughlin's Module}


We recall the definition of the module $M$ given by Dutta, Hochster and
McLaughlin \cite{dut-hoc-mcl:mfpdnim}. As a $\kappa$-vector space, it is
$15$-dimensional:
$$
M=(\kappa u_1+\cdots+\kappa u_{5})+(\kappa v_1+\cdots+\kappa v_{4})+
(\kappa w_1+\cdots+\kappa w_{6}).
$$
Its module structure is given by
\begin{eqnarray*}
Xu_i=Yu_i=Zu_i=Wu_i=0 & & (i=1,\cdots,5)
\end{eqnarray*}
\begin{equation*}
\begin{array}{lllllll}
Xv_1   = u_1 & & Yv_1  = 0     & & Zv_1  = 0             & & Wv_1   = 0 \\
Xv_2   = u_2 & & Yv_2  = 0     & & Zv_2  = 0             & & Wv_2   = 0 \\
Xv_3   = 0     & & Yv_3  = 0     & & Zv_3  = u_1         & & Wv_3   = 0 \\
Xv_4   = 0     & & Yv_4  = 0     & & Zv_4  = u_2         & & Wv_4   = 0
\end{array}
\end{equation*}
\begin{equation*}
\begin{array}{lllllll}
Xw_1 = v_1 & & Yw_{1} = u_3 & & Zw_{1} = 0             & & Ww_{1} = u_1 \\
Xw_2 = v_2 & & Yw_{2} = u_4 & & Zw_{2} = 0             & & Ww_{2} = u_2 \\
Xw_3 = v_3 & & Yw_{3} = u_5 & & Zw_{3} = v_1         & & Ww_{3} = 0     \\
Xw_4 = v_4 & & Yw_{4} = 0     & & Zw_{4} = v_2         & & Ww_{4} = u_3 \\
Xw_5 = u_4 & & Yw_{5} = 0     & & Zw_{5} = v_3         & & Ww_{5} = u_4 \\
Xw_6 = u_5 & & Yw_{6} = 0     & & Zw_{6} = u_3+v_4 & & Ww_{6} = u_5.
\end{array}
\end{equation*}
Note that all monomials of degree greater than one act on the basis 
$u_i$, $v_j$, $w_k$ trivially
except the following cases.
$$
\begin{array}{lllll}
X^2w_1=u_1  & & XZw_3 = u_1  &  & Z^2w_5=u_1 \\
X^2w_2=u_2  & & XZw_4 = u_2  &  & Z^2w_6=u_2
\end{array}
$$

An $A$-linear homomorphism $\Phi$ from $M$ to an $A$-module $N$ is 
determined by its
values at $w_1,\cdots w_6$ and satisfies the conditions
$$
Z\Phi(w_1)=Z\Phi(w_2)=W\Phi(w_3)=Y\Phi(w_4)=Y\Phi(w_5)=Y\Phi(w_6)=0
$$
$$
\begin{array}{rllllll}
X\Phi(w_1) &=& Z\Phi(w_3) \\
X\Phi(w_2) &=& Z\Phi(w_4) \\
X\Phi(w_3) &=& Z\Phi(w_5) \\
& &Y\Phi(w_1) &=& W\Phi(w_4) \\
X\Phi(w_5) &=& Y\Phi(w_2) &=& W\Phi(w_5)\\
X\Phi(w_6) &=& Y\Phi(w_3) &=& W\Phi(w_6)
\end{array}
$$
$$
Z\Phi(w_6)=Y\Phi(w_1)+X\Phi(w_4)
$$
\begin{eqnarray*}
W\Phi(w_1)=X^2\Phi(w_1)=XZ\Phi(w_3)=Z^2\Phi(w_5)\\
W\Phi(w_2)=X^2\Phi(w_2)=XZ\Phi(w_4)=Z^2\Phi(w_6)
\end{eqnarray*}
and the condition that all monomials of degree greater than one act 
trivially on $\Phi(w_i)$
except $X^2\Phi(w_1)$, $XZ\Phi(w_3)$, $Z^2\Phi(w_5)$, $X^2\Phi(w_2)$,
$XZ\varphi(w_4)$, $Z^2\varphi(w_6)$.
Any six elements $\Phi(w_1),\cdots,\Phi(w_6)\in N$ satisfying the above
conditions extend uniquely to an $A$-linear map $\Phi\colon M\to N$.
Note that some of these conditions are redundant.
\begin{lem}\label{lem:573928}
$\Hom_A(M,E(f))=0$.
\end{lem}
\begin{proof}
Let $\Phi\in\Hom_A(M,E(f))$. If $Z\not\in(f)$, multiplication by $Z$ 
is bijective.
Thus
\begin{eqnarray*}
Z\Phi(\omega_1)=Z\Phi(\omega_2)=0&\implies&\Phi(\omega_1)=\Phi(\omega_2)=0\\
Z\Phi(\omega_3)=X\Phi(\omega_1)=0&\implies&\Phi(\omega_3)=0\\
Z\Phi(\omega_4)=X\Phi(\omega_2)=0&\implies&\Phi(\omega_4)=0\\
Z\Phi(\omega_5)=X\Phi(\omega_3)=0&\implies&\Phi(\omega_5)=0\\
Z^2\Phi(\omega_6)=W\Phi(\omega_2)=0&\implies&\Phi(\omega_6)=0
\end{eqnarray*}
If $f=Z$, multiplication by $W$ is bijective. Thus
\begin{eqnarray*}
W\Phi(\omega_3)=0&\implies&\Phi(\omega_3)=0\\
W\Phi(\omega_6)=Y\Phi(\omega_3)=0&\implies&\Phi(\omega_6)=0\\
W\Phi(\omega_2)=Z^2\Phi(\omega_6)=0&\implies&\Phi(\omega_2)=0\\
W\Phi(\omega_5)=Y\Phi(\omega_2)=0&\implies&\Phi(\omega_5)=0\\
W\Phi(\omega_1)=Z^2\Phi(\omega_5)=0&\implies&\Phi(\omega_1)=0\\
W\Phi(\omega_4)=Y\Phi(\omega_1)=0&\implies&\Phi(\omega_4)=0
\end{eqnarray*}
In either case, $\Phi=0$.
\end{proof}

Now we compute $M':=\Hom_A(M,E(Z,W))$. For $1\leq i,j \leq 6$, let
$$
\Phi_i(w_j):=\delta_{ij}\Omega^0(1).
$$
Furthermore, for $1\leq j\leq 6$, let
\begin{eqnarray*}
\Phi_{13}(\omega_j)&:=&X^{-1}\Phi_1(\omega_j)+Z^{-1}\Phi_3(\omega_j),\\
\Phi_{24}(\omega_j)&:=&X^{-1}\Phi_2(\omega_j)+Z^{-1}\Phi_4(\omega_j),\\
\Phi_{35}(\omega_j)&:=&X^{-1}\Phi_3(\omega_j)+Z^{-1}\Phi_5(\omega_j),\\
\Phi_{46}(\omega_j)&:=&X^{-1}\Phi_4(\omega_j)+Z^{-1}\Phi_6(\omega_j),\\
\Phi_{14}(\omega_j)&:=&Y^{-1}\Phi_1(\omega_j)+(W^{-1}-X^{-1})\Phi_4(\omega_j),\\
\Phi_{25}(\omega_j)&:=&Y^{-1}\Phi_2(\omega_j)+(W^{-1}+X^{-1})\Phi_5(\omega_j),\\
\Phi_{36}(\omega_j)&:=&Y^{-1}\Phi_3(\omega_j)+(W^{-1}+X^{-1})\Phi_6(\omega_j),\\
\Phi_{135}(\omega_j)&:=&(W^{-1}+X^{-2})\Phi_1(\omega_j)+Z^{-1}X^{-1}\Phi_3(\omega_j)+Z^{-2}\Phi_5(\omega_j),\\
\Phi_{246}(\omega_j)&:=&(W^{-1}+X^{-2})\Phi_2(\omega_j)+Z^{-1}X^{-1}\Phi_4(\omega_j)+Z^{-2}\Phi_6(\omega_j).
\end{eqnarray*}
For $1\leq i \leq 6$, it is straightforward to check that 
$\Phi_i(w_j)$ satisfy the
conditions prior to Lemma~\ref{lem:573928}. For $i \in
\{13,24,35,46,14,25,36,135,246\}$, we use the following facts to
check these conditions for $\Phi_i(w_j)$.
\begin{itemize}
\item Divisions or multiplications by powers of different variables on
   $E(Z,W)$ are commutative. For example, $X^iY^j\Psi=Y^jX^i\Psi$ for
   $i,j\in\mathbb Z$.
\item The divisions by a power of a single variable satisfy the properties:
   $X^iX^{-j}\Psi=X^{i-j}\Psi$ for $i,j>0$ and the same for
   $Y$, $Z$ and $W$. However,
   $X^{-1}X\Psi\neq\Psi$ in general.
\end{itemize}
All these $\Phi_i$ extend to well-defined elements in $M'$.

\begin{lem}
$$
\{\Phi_i|i\in\{1,2,3,4,5,6,13,24,35,46,14,25,36,135,246 \}\}
$$
is a basis for the $\kappa$-vector space $M'$.
\end{lem}
\begin{proof}
Clearly, $\Phi_1$, $\Phi_2$, $\Phi_3$, $\Phi_4$, $\Phi_5$, $\Phi_6$ 
are linearly independent. Assume that
\begin{eqnarray*}
\Phi=\sum_i a_i\Phi_i=0& &(a_i\in\kappa).
\end{eqnarray*}
Since $Z^2\Phi=W\Phi=\Phi=0$, we evaluate the left hand side of the above
equality at $\omega_3$, $\omega_4$, $\omega_5$, $\omega_6$ and obtain
\begin{eqnarray*}\begin{array}{l}
a_{135}=a_{246}=0,\\
a_{14}=a_{25}=a_{36}=0,\\
a_{13}=a_{24}=a_{35}=a_{46}=0.
\end{array}
\end{eqnarray*}
Hence all $a_i=0$ and $\Phi_i$ are linearly independent. Since the 
Matlis dual $M'$ of $M$ has length  $15$,
$\Phi_i$ generate $M'$.
\end{proof}
\begin{prop}
As an $A$-module, the minimal number of generators for $M'$ is $5$.
\end{prop}
\begin{proof}
$(X,Y,Z,W)M'$ as a vector space is generated by $\Phi_1$, $\Phi_2$, $\Phi_3$,
$\Phi_4$, $\Phi_5$, $\Phi_6$, $\Phi_{13}$, $\Phi_{24}$, $\Phi_{35}$, 
$\Phi_{46}$. Therefore
$\{\Phi_{14},\Phi_{25},\Phi_{36},\Phi_{135},\Phi_{246}\}$ is a 
minimal generating set for $M'$.
\end{proof}

Now we compute $\Ext^i_A(M,A/\p)$.

\begin{prop}
$$ \dim_k \Ext^i_A(M,A/\p)= \begin{cases}
6, & \mbox{ if } i=2;\\
7, & \mbox{ if } i=3;\\
0, &\text{otherwise}.
\end{cases}
$$
\end{prop}

\begin{proof}
For $n\geq 2$, $\Ext^{2n}_A(M,A/\p)$ is the cohomology of
$$
{M'}^2\xrightarrow{\tiny\begin{pmatrix}W&-Z\\-Y&X\end{pmatrix}}
{M'}^2\xrightarrow{\tiny\begin{pmatrix}X&Z\\Y&W\end{pmatrix}}{M'}^2.
$$
To simplify the notation, we write $(i,j)$ for the element 
$(\Phi_i,\Phi_j)\in{M'^2}$. The
kernel of $\tiny\begin{pmatrix}X&Z\\Y&W\end{pmatrix}$ is generated by 
$(35,-13)$,
$(46,-24)$, and $(i,0)$, $(0,i)$, where $1\leq i\leq 6$. The image of
$\tiny\begin{pmatrix}W&-Z\\-Y&X\end{pmatrix}$ is generated by 
$(4,-1)$, $(5,-2)$, $(6,-3)$, $(1,0)$,
$(2,0)$, $(-3,1)$, $(-4,2)$, $(-5,3)$, $(-6,4)$, $(0,-4)$, $(0,5)$, 
$(0,6)$, $(-35,13)$, $(-46,24)$. Since
\begin{eqnarray*}
(0,-3)&=&(6,-3)+(-6,4)+(0,-4),\\
(0,-2)&=&(5,-2)+(-5,3)+(0,-3),\\
(0,-1)&=&(4,-1)+(-4,2)+(0,-2),\\
(-3,0)&=&(-3,1)+(0,-1),\\
(-4,0)&=&(-4,2)+(0,-2),\\
(-5,0)&=&(-5,3)+(0,-3),\\
(-6,0)&=&(-6,4)+(0,-4),\\
\end{eqnarray*}
all $(0,i)$ and $(i,0)$ are contained in the image of 
$\tiny\begin{pmatrix}W&-Z\\-Y&X\end{pmatrix}$. Therefore
$\Ext^{2n}_A(M,A/\p)=0$ for $n\geq 2$.

For $n\geq 2$, $\Ext^{2n+1}_A(M,A/\p)$ is the cohomology of
$$
{M'}^2\xrightarrow{\tiny\begin{pmatrix}X&Z\\Y&W\end{pmatrix}}
{M'}^2\xrightarrow{\tiny\begin{pmatrix}W&-Z\\-Y&X\end{pmatrix}}{M'}^2.
$$
The kernel of $\tiny\begin{pmatrix}W&-Z\\-Y&X\end{pmatrix}$ is generated by
$(13,0)$, $(24,0)$, $(35,0)$, $(46,0)$ and $(i,0)$, $(0,i)$, where 
$1\leq i\leq 6$. The image of
$\tiny\begin{pmatrix}X&Z\\Y&W\end{pmatrix}$ is generated by $(1,0)$, 
$(2,0)$, $(3,0)$, $(4,0)$, $(-4,1)$, $(5,2)$,
$(6,3)$, $(13,0)$, $(24,0)$, $(5,0)$, $(6,0)$, $(0,4)$, $(0,5)$, 
$(0,6)$, $(35,1)$, $(46,2)$. Clearly
$\Ext^{2n+1}_A(M,A/\p)=0$ for $n\geq 2$.

$\Ext^3_A(M,A/\p)$ is the cohomology of
$$
M'\xrightarrow{\tiny\begin{pmatrix}X\\Y\end{pmatrix}}
{M'}^2\xrightarrow{\tiny\begin{pmatrix}W&-Z\\-Y&X\end{pmatrix}}{M'}^2.
$$
The image of $\tiny\begin{pmatrix}X\\Y\end{pmatrix}$ is generated by 
$(1,0)$, $(2,0)$, $(3,0)$, $(4,0)$, $(-4,1)$, $(5,2)$,
$(6,3)$, $(13,0)$, $(24,0)$. $\Ext^3_A(M,A/\p)$ is generated by the 
classes represented by
$(0,2)$, $(0,3)$,
$(0,4)$, $(0,5)$, $(0,6)$, $(35,0)$, $(46,0)$. In particular,
$\dim_\kappa\Ext^3_A(M,A/\p)=7$.

$\Ext^2_A(M,A/\p)$ has a basis $\Phi_1,\cdots,\Phi_6$,
so $\dim_k \Ext^2_A(M,A/\p) =6$. It is easy to see that
$\Ext^0_A(M,A/\p)=\Ext^1_A(M,A/\p)=0$. This completes the proof of the
proposition.

\end{proof}


\end{document}